\documentclass[11pt]{article}

\usepackage[margin=1in]{geometry}
\usepackage{amsmath,amssymb,amsthm}
\usepackage{booktabs}
\usepackage{array}        
\usepackage{graphicx}
\usepackage{microtype}
\usepackage[numbers,sort&compress]{natbib}
\usepackage[hidelinks]{hyperref}
\usepackage{xurl}  
\usepackage{xcolor}
\usepackage[font=small,labelfont=bf]{caption}
\usepackage{subcaption}
\usepackage{float}

\graphicspath{{figures/}}

\definecolor{doilink}{RGB}{0,76,153}
\newcommand{\doilink}[1]{\href{https://doi.org/#1}{\textcolor{doilink}{\texttt{#1}}}}
\newcommand{\VE}{\mathrm{VE}}
\newcommand{\dist}{D}
\usepackage{authblk} 

\theoremstyle{definition}
\newtheorem{proposition}{Proposition}
\theoremstyle{definition}
\newtheorem{definition}{Definition}

\newcommand{\Gr}{\mathrm{Gr}}
\newcommand{\St}{\mathrm{St}}
\newcommand{\alphabet}{\mathcal{A}}
\newcommand{\chan}{\mathcal{G}}
\newcommand{\rank}{\mathrm{rank}}
\newcommand{\tr}{\mathrm{tr}}
\newcommand{\dGr}{d_{\Gr}}
\newcommand{\dbeta}{d_{\beta}}

\title{Mathematical AI Predicts 40\% and 57\% Influenza Vaccine Effectiveness Against H1N1 and H3N2, Respectively, for the 2026–2027 U.S. Season}

\author[1]{ Xiang Xiang Wang }
\author[1,2,3]{Guo-Wei Wei\thanks{Corresponding author: Guo-Wei Wei (\texttt{guowei.wei@uga.edu}).}}

\affil[1]{Department of Mathematics, University of Georgia, Athens, GA 30602, USA}
\affil[2]{Department of Biochemistry and Molecular Biology, University of Georgia, Athens, GA 30602, USA}
\affil[3]{School of Computing, University of Georgia, Athens, GA 30602, USA}

\date{}

\begin{document}
\maketitle

\begin{abstract}
\noindent
\textbf{Background.}
Clinical influenza vaccine effectiveness (VE) has often been low over the past two decades, highlighting the need for preseason VE forecasting. Timely forecasts could  
help health authorities plan prevention strategies and allocate healthcare resources.
 However, preseason forecasting of clinical influenza VE remains relatively uncommon.
 
\textbf{Methods.}
We obtained hemagglutinin (HA) sequences from GISAID for U.S. specimens collected in August--September before each influenza season from 2010--2011 through 2024--2025. After quality screening, we compared each virus with the corresponding vaccine target strain for the same season and subtype or lineage using two mathematical distances based on amino acid patterns in the HA protein. We calculated mean distances for each season and subtype or lineage and paired them with published clinical VE estimates, yielding 26 observations across 13 seasons. Using a mathematical AI approach, we fitted linear regression models of VE against each distance and evaluated their predictive performance by forecasting later seasons using data from earlier seasons, with the historical mean VE for each subtype or lineage as a reference. We then applied the same approach to 87 HA sequences collected in August--September 2026 to forecast VE for the 2026--2027 season.

\textbf{Results.}
We generated component-specific forecasts of influenza VE for the upcoming 2026--2027 U.S. season. Both mathematical distances were negatively associated with historical VE, indicating that greater differences between circulating viruses and vaccine strains were associated with lower VE. The topological model forecast VE of 40.3\%, 56.9\%, and 47.1\% against A(H1N1)pdm09, A(H3N2), and B/Victoria, respectively, whereas the geometric model forecast corresponding VE values of 34.0\%, 56.8\%, and 51.1\%. 


\end{abstract}

\section{Introduction}

Influenza   vaccination is an important means of reducing seasonal illness, but
its effectiveness varies from one season to another \citep{Belongia2016}. Knowing the likely level of protection before the influenza season would help public health authorities prepare for a potentially substantial influenza burden, plan healthcare capacity, refine risk communication, and inform public health decision-making, particularly when protection against a specific subtype is expected to be limited. Clinical vaccine effectiveness (VE)
estimates, however, require vaccination records and laboratory test results
from people who develop respiratory illness and seek medical care
\citep{JacksonNelson2013}. These data
accumulate as the season progresses. Interim estimates can be reported during
the season: the 2025--2026 US interim report was published on 12 March 2026
using surveillance data collected from September 2025 through February 2026
\citep{Maloney2026}. Full-season journal reports may appear substantially
later. For example, reports for the 2017--2018, 2023--2024 and 2024--2025
seasons were first published online on 2 February 2019, 6 January 2025 and
14 July 2026, respectively \citep{Rolfes2019,Chung2025,Chung2026}. These dates
refer to journal publication rather than the first public availability of any
VE information; the  Centers for Disease Control and Prevention (CDC) also presented preliminary end-of-season estimates for
2025--2026 in June 2026 \citep{Grohskopf2026Presentation}. These estimates
measure protection using observed clinical outcomes rather than forecast
protection before those outcomes occur. Earlier forecasts could complement clinical
surveillance by providing additional time to prepare for lower protection,
with the potential to reduce infections and avoidable healthcare costs.

Virus sequences offer an information source that can be examined before sufficient
clinical outcomes have accumulated. Influenza viruses undergo genetic and
antigenic evolution \citep{Smith2004}, and changes at specific HA positions can
alter antigenic properties \citep{Koel2013}. Differences between circulating
and vaccine viruses are relevant to protection: the 2025--2026 US interim
report documented the predominance of antigenically drifted A(H3N2) viruses
alongside lower interim VE than in recent seasons \citep{Maloney2026}. Hemagglutinin (HA) sequences describe variation in a major
surface protein. Experimental studies have linked changes at specific HA
positions to antigenic changes \citep{Koel2013}, and sequence-based antigenic
distances have been related to historical vaccine protection \citep{Gupta2006}.  However, larger sequence differences do not necessarily imply larger antigenic differences \citep{Smith2004}. VE also depends on factors such as prior immunity,
vaccine characteristics, the study population and the clinical outcome
assessed \citep{Chung2025,Chung2026}. A key challenge is determining whether distances between preseason HA
sequences and vaccine strains can provide useful information for predicting VE.

Mathematical artificial intelligence (AI), particularly  topological deep learning  \citep{CangWei2017}, offers an effective  approach for biomolecular analysis \citep{WeeWei2025}.
Mathematical representations of biomolecular data of this kind are reviewed by
\citet{NguyenCangWei2020}, and persistent topological Laplacians have been
applied to the structural analysis of viral proteins \citep{WeiChenWei2023}. 
Topological deep learning and spectral representations have also been used to study how viral
mutations affect viral infectivity and antibody resistance \citep{ChenWangWangWei2020,WangChenWei2021,ChenGilbyWei2022}. These approaches lead to the successful predictions of crucial spike protein mutation sites  \citep{ChenWangWangWei2020} and emerging SARS-CoV-2 dominant  variants 
\citep{chen2022Omicron, ChenQiuWangWei2022}. 
More recently, Hozumi and Wei introduced $k$-mer topology to compare sequences
through persistent topological and spectral features, and examined its use for
antigenic-distance analysis of SARS-CoV-2 variants \citep{HozumiWei2024}.

In this work, we propose the use of mathematical distances between vaccine target and influenza variants to predict influenza VE. Our study combines HA sequences from GISAID  \citep{ShuMcCauley2017}  with published clinical VE estimates to forecast VE for a future season. We use topological and geometric methods to quantify differences between
preseason HA sequences and vaccine strains, then relate these distances to
historical clinical VE through regression. 
We construct two distances from the
positions of short amino-acid fragments, or $k$-mers, along HA. The topological
distance compares how occurrences of each fragment form connected groups
across scales \citep{HozumiWei2024}. The geometric distance represents
topological and Laplacian spectral descriptors as subspaces on a Grassmann
manifold and compares these subspaces using chordal distance
\citep{WangNguyenWei2020,WangWei2025}. Grassmann-manifold representations of
this kind have also been used for other high-dimensional biological data
\citep{WangCottrellWei2026}, and related persistent-Laplacian constructions
extend to sheaf-theoretic settings \citep{WangWeiSheaf2026}. These constructions summarize sequence
variation in different ways. Evaluating them on the same sequences, vaccine
targets and clinical observations allows us to compare their relationships
with VE.

Using August--September sequences to forecast next season VE requires addressing two
questions. First, do virus--vaccine target distances measured before the season
reflect those observed during the subsequent clinical enrollment period?
Second, can the relationship between sequence distance and historical VE
help predict effectiveness in a later season? These questions are related
but distinct. Agreement between collection windows supports the timing of
the sequence measurement, whereas prediction against clinical outcomes
assesses its value for forecasting protection. Examining both provides a
basis for evaluating preseason sequence distances as early indicators of VE.

We pair US HA sequences from specimens collected during August and September
with published component-specific VE estimates of the next season, yielding 26 season--component
observations from 13 historical influenza seasons between 2010--2011 and
2024--2025. After screening sequences and
verifying the clinical estimates, we calculate distances to the corresponding
vaccine strains and fit a separate regression for each distance. We assess prediction on later seasons using pooled models and distance
normalisation fitted to earlier data. We also compare topological distances calculated from
August--September sequences with those from sequences collected during the
subsequent VE-study enrollment window, to assess whether preseason measurements
reflect later-season virus--vaccine differences. We then forecast VE against
A(H1N1)pdm09, A(H3N2) and B/Victoria
for the 2026--2027 US influenza season using 87 eligible HA sequences from
specimens collected in the United States during August--September 2026.
These forecasts are intended for comparison with final US Flu VE Network
estimates for medically attended, laboratory-confirmed outpatient influenza,
matched by subtype or lineage. This comparison target differs from the
aggregate overall VE on CDC's seasonal summary webpage.

\section{Mathematical foundations}
\label{sec:framework}

We describe each HA sequence in two steps. First, the positions at which
short fragments recur define a family of graphs: increasing the allowed gap
joins separate occurrences into connected groups. Second, the resulting
measurements are compared either as vectors or as linear subspaces. The
first view retains differences in their numerical values; the second
compares the subspaces they span. We first introduce the topological objects, then construct them from
$k$-mer positions, before describing the geometry of subspaces and the two
distance functions.

\subsection{Basic topological concepts}
\label{sec:topologybasics}

\begin{figure}[!htb]
  \centering
  \captionsetup{skip=4pt}
  \includegraphics[width=\textwidth]{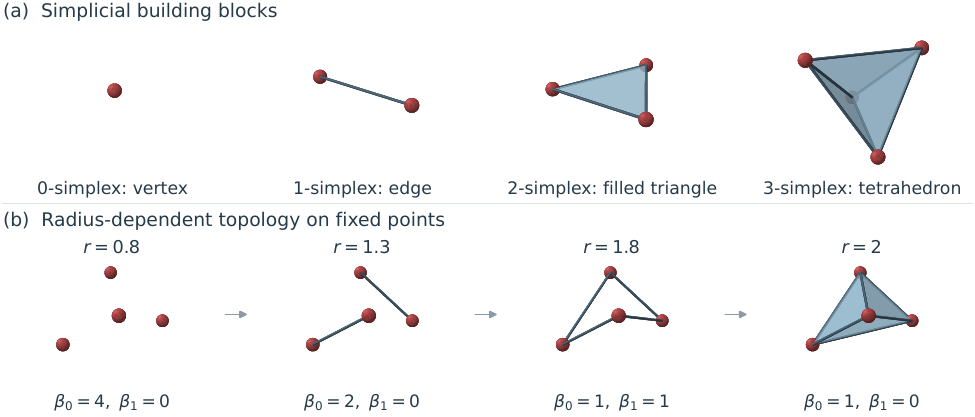}
  \caption{\textbf{Simplices and a radius-based filtration.}
  (a) Basic simplices of dimensions zero to three.
  (b) Vietoris--Rips complexes on the same four non-coplanar points,
  shown from the same three-dimensional viewpoint. An edge enters when its
  Euclidean length is at most $r$, and every clique defines a simplex.
  Increasing $r$ joins four isolated points into two components and then
  one loop; adding faces fills the loop, and the final complex is a solid
  tetrahedron. The displayed Betti numbers count components and independent
  loops. Shaded faces show the three-dimensional structure; sphere and edge
  thicknesses are used only for visibility. All panels schematically illustrate basic topological concepts.}
  \label{fig:topologybasics}
\end{figure}

Topology describes features such as connectivity and holes. For finite data,
simplicial complexes provide a way to describe these features using simple
building blocks \citep[Chapter III]{EdelsbrunnerHarer2010}. A vertex is a 0-simplex, an edge is a 1-simplex, a filled
triangle is a 2-simplex, and a solid tetrahedron is a 3-simplex
(Figure~\ref{fig:topologybasics}a). These standard examples also introduce the
basic objects in the persistent-Laplacian treatment of \citet{WeiWei2021}.

More generally, a $q$-simplex is the convex hull of $q+1$ affinely
independent vertices $v_0,\ldots,v_q$:
\[
 [v_0,\ldots,v_q]
 =\left\{\sum_{i=0}^{q}t_i v_i:
          t_i\geq0,\ \sum_{i=0}^{q}t_i=1\right\}.
\]
A face is the simplex obtained from a subset of these vertices. A simplicial
complex is a collection of simplices containing all their faces, with any
nonempty intersection shared as a common face. Its 1-skeleton consists of
its vertices and edges and is a graph. Vertices belong to the same connected
component when an edge path joins them.

Homology describes connectivity and higher-dimensional holes algebraically
\citep[Chapter IV]{EdelsbrunnerHarer2010}.
Over a field $\mathbb{F}$, the $q$th Betti number is
$\beta_q=\dim H_q(\mathcal{K};\mathbb{F})$ for a complex $\mathcal{K}$.
In low dimensions, $\beta_0$ counts connected components, $\beta_1$ counts
independent loops and $\beta_2$ counts enclosed cavities. An unfilled
triangular boundary has one loop; adding its triangular face fills that
loop. A simplex itself has one connected component and no higher-dimensional
holes. These distinctions explain why vertices, edges and filled faces
must be specified separately.

The scale at which points are connected determines which of these features
are visible. To examine more than one scale, let
$X=\{x_1,\ldots,x_n\}\subset\mathbb{R}^d$ be a finite point cloud.
Its Vietoris--Rips complex connects points according to a distance threshold
\citep[Section III.2]{EdelsbrunnerHarer2010}. Here, $r$ denotes the edge-length
threshold, so the complex contains precisely the vertex sets whose pairwise
distances are at most $r$:
\[
 \operatorname{VR}_r(X)
 =\{\sigma\subseteq X:\|x-y\|\leq r
                    \text{ for every }x,y\in\sigma\}.
\]
Thus a close pair supplies an edge, three mutually close points supply a
filled triangle, and four mutually close points supply a 3-simplex. This
rule specifies the faces as well as the edges; drawing an edge cycle alone
is not the same as filling its interior.

For $r_0\leq r_1\leq\cdots\leq r_m$, every simplex present at one radius
remains present at larger radii, giving a filtration
\citep[Chapter VII]{EdelsbrunnerHarer2010}
\[
 \operatorname{VR}_{r_0}(X)\subseteq
 \operatorname{VR}_{r_1}(X)\subseteq\cdots\subseteq
 \operatorname{VR}_{r_m}(X).
\]
The point coordinates remain fixed throughout; it is their connection
structure that changes. Figure~\ref{fig:topologybasics}b illustrates this
process on four points. At $r=0.8$ they form four separate components.
At $r=1.3$, two edges leave two components. At $r=1.8$, a connected cycle
has formed, so $\beta_0=1$ and $\beta_1=1$. At $r=2$, the additional edges
and faces fill the loop and complete a tetrahedron, leaving $\beta_0=1$
and $\beta_1=0$. This is the radius-dependent construction illustrated in
the persistent-Laplacian examples of \citet{WeiWei2021}.

Persistent homology tracks the appearance and disappearance of homological
features along this filtration \citep[Chapter VII]{EdelsbrunnerHarer2010}. A feature is born when it first appears and
dies when it merges with an older component or becomes a boundary. For
example, adding edges can merge components, whereas adding triangular
faces can remove loops. The interval between birth and death records
persistence across connection scales, rather than elapsed biological time.
Recording only $\beta_q(r)$ gives a Betti curve: a count of the features
present at each radius, without retaining their individual identities.

These examples introduce the general construction. In the present study,
the points are occurrence positions along an HA sequence, and the distances
use component counts and graph-Laplacian measurements across radii. The
loops and filled faces above explain the broader topological language;
higher-dimensional Betti numbers are not inputs to the VE models. The next
subsection specifies the one-dimensional position sets and their graphs.

\subsection{$k$-mer positions and their topology}
\label{sec:pointclouds}

Following the position-based representation of \citet{HozumiWei2024}, let $\alphabet$ be a finite alphabet with $A = |\alphabet|$, and let
$S = S_1 S_2 \cdots S_L \in \alphabet^{L}$ be a sequence. For an integer
$k \geq 1$ the set of $k$-mer types is $\alphabet^{k}$, of cardinality
$A^{k}$. For a $k$-mer type $l \in \alphabet^{k}$ define its occurrence set
\begin{equation}
  P^{l}(S) \;=\; \bigl\{\, i \in \{1,\dots,L-k+1\} \;:\;
  S_i S_{i+1} \cdots S_{i+k-1} = l \,\bigr\} \;\subset\; \mathbb{Z},
  \label{eq:occurrences}
\end{equation}
a finite point cloud on the integer line recording where in $S$ that
$k$-mer occurs. A window containing a character outside $\alphabet$ --- a
nucleotide ambiguity code, or a non-standard residue --- matches no type and so
belongs to no occurrence set. The indexed collection $\{P^{l}(S)\}_{l \in \alphabet^{k}}$
records the positions of each $k$-mer type without requiring an alignment
between sequences.

This construction applies to any finite alphabet. Here, $\alphabet$ is
the set of 20 standard amino acids; an individual amino acid, an adjacent
pair and an adjacent triplet correspond to $k=1,2,3$. The source methods
provide the position-based construction and its nucleotide implementations
\citep{HozumiWei2024,WangWei2025}; the protein settings used here are specified
in Section~\ref{sec:instantiation}.

\label{sec:betti}

For a radius $r\geq0$, join two distinct positions $i,j\in P^l(S)$ when
$|i-j|\leq r$. This gives a graph $G^l_r(S)$, the 1-skeleton of the
Vietoris--Rips complex at scale $r$. Positions are compared along the sequence
as written; \citet{HozumiWei2024} also describe a circular variant that joins
the two ends, which is not used here. Every set of pairwise connected vertices specifies a simplex, but the
vertices and edges alone determine the component count. Increasing $r$
produces the filtration described above: $G^l_r(S)\subseteq G^l_{r'}(S)$
whenever $r\leq r'$.

The zeroth homology space $H_0(G^l_r(S);\mathbb{F})$, over a field
$\mathbb{F}$, has one independent class for each component. Its dimension
is therefore the component count. Equivalently, let
$L^l_r(S)=\mathsf{D}-\mathsf{A}$ be the graph Laplacian, where
$\mathsf{A}$ is the adjacency matrix and $\mathsf{D}$ contains the vertex
degrees on its diagonal \citep{Merris1994}. The multiplicity of the zero
eigenvalue equals the number of connected components, giving
\begin{equation}
  \beta_0^{l}(S;r)
  =\dim H_0(G^l_r(S);\mathbb{F})
  =\#\{\text{connected components of }G^l_r(S)\}
  =\dim\ker L^l_r(S).
  \label{eq:betti0}
\end{equation}
The function $r\mapsto\beta_0^l(S;r)$ is the Betti-0 curve. For a nonempty
occurrence set it decreases from the number of occurrences to one. For positions $\{1,3,8\}$, the gap of two joins the first two points,
and the gap of five joins the remaining point. The curve thus retains
information about separation as well as occurrence count.

Because the point cloud lies in one dimension, this curve has a closed form.

\begin{samepage}
\begin{proposition}[Betti-0 of a one-dimensional point cloud]
\label{prop:betti}
Let $p_1 < p_2 < \cdots < p_n$ be the elements of $P^{l}(S)$ in increasing
order and let $g_i = p_{i+1} - p_i$ for $1 \leq i \leq n-1$ be the consecutive
gaps. Then for every $r \geq 0$,
\begin{equation}
  \beta_0^{l}(S;r) \;=\; 1 + \#\{\, i : g_i > r \,\}
  \qquad (n \geq 1),
  \label{eq:betticlosed}
\end{equation}
and $\beta_0^{l}(S;r) = 0$ when $n = 0$.
\end{proposition}

The position-based filtration is that of \citet{HozumiWei2024}; in one
dimension each consecutive gap larger than $r$ separates two components, which
gives the closed form above.
\end{samepage}

Proposition~\ref{prop:betti} allows the Betti-0 curve to be evaluated
without constructing a simplicial complex or a graph Laplacian. After sorting
the positions and their consecutive gaps in $O(n \log n)$ time, evaluation
at $T$ radii costs $O(T \log n)$ by binary search on the sorted gaps.

\label{sec:laplacian}

The Betti-0 curve reads only the harmonic part of the spectrum --- the
multiplicity of the eigenvalue $0$. The non-harmonic spectrum, consisting of the positive eigenvalues
\citep{WangNguyenWei2020,WeiWei2021}, carries complementary information about how tightly each cluster
of occurrences is connected, not merely how many clusters there are. Write
$m = n - \beta_0^{l}(S;r) = \rank L^{l}_{r}(S)$ for the number of strictly
positive eigenvalues and, when $m \geq 1$, define the first two moments of the
non-harmonic spectrum
\begin{equation}
  \mu^{l}(S;r) = \frac{1}{m}\!\!\sum_{\lambda > 0}\!\! \lambda,
  \qquad
  \sigma^{l}(S;r) = \Bigl(\tfrac{1}{m}\!\!\sum_{\lambda > 0}\!\! \lambda^{2}
  - \mu^{l}(S;r)^{2}\Bigr)^{1/2},
  \label{eq:moments}
\end{equation}
the sums running over the positive eigenvalues of $L^{l}_{r}(S)$ with
multiplicity. These too avoid an eigendecomposition.

\begin{samepage}
\begin{proposition}[Degree-only form of the spectral moments]
\label{prop:moments}
Let $\deg_i$ be the degree of vertex $i$ in $G^{l}_{r}(S)$. Then
\begin{equation}
  \sum_{\lambda > 0} \lambda \;=\; \tr L^{l}_{r}(S) \;=\; \sum_{i} \deg_i,
  \qquad
  \sum_{\lambda > 0} \lambda^{2} \;=\; \tr \bigl(L^{l}_{r}(S)\bigr)^{2}
  \;=\; \sum_{i} \deg_i^{2} + \sum_{i} \deg_i .
  \label{eq:traces}
\end{equation}
Consequently $\mu^{l}$ and $\sigma^{l}$ are determined by the degree sequence
alone and, after sorting, are computable in $O(n)$ per radius.
\end{proposition}

These identities follow from the standard trace relations for the combinatorial graph Laplacian \citep{Merris1994}.
\end{samepage}

The smallest positive eigenvalue $\lambda_{\min}^{+}(L^l_r(S))$ supplies
a further spectral summary \citep{Merris1994,WangNguyenWei2020}. When several nontrivial components are present,
it is the minimum of their smallest positive eigenvalues. Implementation
details are given in Supplementary Information, Section~\ref{sec:implementation}.

\label{sec:blocks}

Fix a finite grid of radii $r_1 < r_2 < \cdots < r_T$ and a set of
channels $\chan$, each channel being one of the scalar summaries above:
\begin{equation}
  \chan \;\subseteq\; \bigl\{\, \beta_0,\;
  \lambda_{\min}^{+},\; \mu,\; \sigma \,\bigr\}.
  \label{eq:channels}
\end{equation}
For a sequence $S$, a channel $g \in \chan$ and an order $k$, define the
feature block
\begin{equation}
  G^{(k)}_{S,g} \in \mathbb{R}^{A^{k} \times T},
  \qquad
  \bigl( G^{(k)}_{S,g} \bigr)_{l,t} \;=\; g\bigl( P^{l}(S);\, r_t \bigr),
  \label{eq:block}
\end{equation}
one row per $k$-mer type in lexicographic order and one column per radius. The
collection $\{G^{(k)}_{S,g}\}_{g \in \chan,\, k \leq K}$ is the common object
from which both predictors are derived. Each matrix has a word-type axis and a radius axis. This common arrangement
lets us compare its entries directly or compare the subspaces obtained after
combining the measurements.

The blocks are sparse: when $A^{k} \gg L$, a $566$-residue protein contains at most $564$ distinct
$3$-mers out of $20^{3} = 8{,}000$ possible types, so most rows vanish
identically. Rows that vanish for every sequence being compared may be omitted from
sparse computations, provided their common coordinates and the row-count
convention used in normalisation are preserved. The entries of different channels and different $k$ have
incommensurable units and magnitudes, so the blocks must be normalised before they are combined; both distances below do so, and both therefore require a scale
estimated over a corpus of comparisons rather than over a single pair.

\subsection{Subspaces and Grassmannian geometry}
\label{sec:grassmann}

A nonzero vector determines a line through the origin; two independent
vectors determine a plane. Changing the lengths or the basis vectors need
not change that line or plane. Grassmannian geometry treats the subspace
itself as the object, separating it from the particular vectors used to
represent it. We use the standard definitions below
\citep{Edelman1998,Conway1996}.

\begin{definition}[Grassmann manifold]
For $1\leq\rho<D$, the real Grassmann manifold is \citep{Edelman1998}
\[
 \Gr(\rho,D)=\{\mathcal{U}\subseteq\mathbb{R}^D:
       \mathcal{U}\text{ is a linear subspace},\ \dim\mathcal{U}=\rho\}.
\]
A point of this manifold is an entire $\rho$-dimensional subspace, rather
than a point lying in that subspace.
\end{definition}

Thus $\Gr(1,3)$ contains the unoriented lines through the origin in
$\mathbb{R}^3$, while $\Gr(2,3)$ contains the planes through the origin.
More generally, $\Gr(\rho,D)$ is a compact smooth manifold of dimension
$\rho(D-\rho)$ \citep{Edelman1998}.

To perform calculations, choose an orthonormal basis for the subspace and
place its vectors in the columns of $U$. Such matrices form the Stiefel
manifold \citep{Edelman1998}
\[
 \St(\rho,D)=\{U\in\mathbb{R}^{D\times\rho}:U^\top U=I_\rho\}.
\]
The columns are mutually perpendicular and have unit length, but this
representation is not unique. Matrices $U$ and $V$ span the same subspace
exactly when $V=UO$ for an orthogonal matrix
$O\in O(\rho)=\{O\in\mathbb{R}^{\rho\times\rho}:O^\top O=I_\rho\}$.
Consequently,
\[
 \Gr(\rho,D)=\St(\rho,D)/O(\rho).
\]
The quotient identifies all orthonormal bases of a common subspace \citep{Edelman1998}.
Rotating an orthonormal basis within its plane does not move the
Grassmann point.

The orthogonal projector gives a basis-independent representation
\citep[Section 2.5]{GolubVanLoan2013}:
\[
 P_{\mathcal{U}}=UU^\top,
 \qquad P_{\mathcal{U}}^\top=P_{\mathcal{U}},
 \qquad P_{\mathcal{U}}^2=P_{\mathcal{U}},
 \qquad \operatorname{rank}P_{\mathcal{U}}=\rho.
\]
For any $x\in\mathbb{R}^D$, $P_{\mathcal{U}}x$ is the closest point in
$\mathcal{U}$ to $x$. Because $(UO)(UO)^\top=UU^\top$, every choice of
orthonormal basis produces the same projector.

Two subspaces of the same dimension can also be compared through their
principal angles \citep[Section 6.4.3]{GolubVanLoan2013}. For orthonormal bases $U,V$, let
$\eta_1\geq\cdots\geq\eta_\rho$ be the singular values of $U^\top V$.
The principal angles are
\[
 \theta_i=\arccos\eta_i,
 \qquad0\leq\theta_1\leq\cdots\leq\theta_\rho\leq\pi/2.
\]
For two lines this becomes the acute angle
$\theta=\arccos|u^\top v|$ between unit representatives; reversing
either vector does not change it.
For higher-dimensional subspaces, successive angle pairs are chosen
orthogonally to the preceding pairs within each subspace. All angles are
zero when the subspaces coincide. These projectors and angles provide two
equivalent ways to express the chordal distance below.

\subsection{Betti-curve Euclidean distance}
\label{sec:reductions}

The first distance compares sequences directly through their Betti-0 curves,
using only the channel $\chan=\{\beta_0\}$. At each radius, a curve records
how many groups of occurrences remain disconnected for a particular
$k$-mer type. At radius zero this is its occurrence count; as the radius
increases, the curve also reflects the gaps separating those occurrences.
For example, positions $\{1,3,8\}$ have consecutive gaps $2$ and $5$, so their
Betti-0 curve equals $3$ for $0\leq r<2$, $2$ for $2\leq r<5$, and $1$ for
$r\geq5$. Comparing these curves therefore incorporates both occurrence
counts and spatial separation along the sequence.

For a fixed order $k$, the block $G^{(k)}_{S,\beta_0}$ contains one row per
$k$-mer type and one column per sampled radius. The operation
$\mathrm{vec}$ arranges all entries into a single vector, with the same
ordering used for every sequence. The within-order distance is
\begin{equation}
  \dbeta^{k}(S,T) \;=\;
  \bigl\lVert \mathrm{vec}\,G^{(k)}_{S,\beta_0}
            - \mathrm{vec}\,G^{(k)}_{T,\beta_0} \bigr\rVert_{2},
  \label{eq:dbetak}
\end{equation}
where the Euclidean norm takes the square root of the sum of squared
entrywise differences. Equivalently, this is the Frobenius norm of the
difference between the two feature blocks \citep[Sections 2.2--2.3]{GolubVanLoan2013}. Each matched $k$-mer--radius
coordinate contributes to the comparison; larger discrepancies contribute
more strongly through squaring.

Different orders capture patterns at different word lengths, and their
raw distances need not have comparable magnitudes. We therefore combine
them as
\begin{equation}
  \dist_{\mathrm{A}}(S,T) \;=\; \sum_{k=1}^{K} a_k \, \dbeta^{k}(S,T),
  \qquad a_k = \frac{1}{\omega_k \, 2^{\,K-k}},
  \label{eq:kmer}
\end{equation}
where $\omega_k>0$ is the mean of $\dbeta^{k}$ over a specified corpus of
comparisons. Division by $\omega_k$ expresses each distance relative to its
reference magnitude. The additional factor $2^{-(K-k)}$ assigns progressively
greater weight to higher orders: it doubles when $k$ increases by one.
These two factors have distinct roles, normalising the distance scale and
setting the relative emphasis across orders. This construction follows the
$k$-mer topological distance of \citet{HozumiWei2024} in their Betti-only
setting, with the spectral term switched off. Their own antigenic-distance
application retains a spectral channel alongside the Betti curves, which is one
reason we also evaluate a distance built from the nonzero spectrum. The scale-fitting populations used in this
study are specified in the methods.

Because the weights are positive and each within-order comparison is a
Euclidean distance, $\dist_{\mathrm A}$ is nonnegative and symmetric. It does
not separate sequences: distinct sequences can have identical sampled Betti
curves, so a zero distance means equality of the chosen topological features
rather than of the sequences themselves. In particular, the curves summarise
connectivity rather than retaining the absolute locations of every occurrence. The construction preserves differences in feature
amplitudes but does not incorporate the nonzero Laplacian spectrum, which
enters the Grassmann representation in the next subsection.

\subsection{Grassmann chordal distance}
To combine $k$-mer topology with Grassmannian geometry, take all four
channels in \eqref{eq:channels}. Normalise each block by a corpus-level root-mean-square
scale $s_{g,k}$ and an order weight $w_k = 2^{-(K-k)}$, and stack the
normalised blocks vertically over $(g,k)$ into one matrix
\begin{equation}
  Z_S \;=\;
  \begin{bmatrix} c_{g,k}\, G^{(k)}_{S,g} \end{bmatrix}_{g \in \chan,\, k \leq K}
  \;\in\; \mathbb{R}^{D \times T},
  \qquad
  c_{g,k} = \frac{\sqrt{w_k}}{s_{g,k}},
  \qquad
  D = |\chan| \sum_{k=1}^{K} A^{k}.
  \label{eq:Z}
\end{equation}
Compute a singular value decomposition of $Z_S$ \citep[Section 2.4]{GolubVanLoan2013} and let
$U_S \in \mathbb{R}^{D \times \rho}$ collect the $\rho$ leading left singular
vectors. The representation of $S$ is the subspace they span,
\begin{equation}
  \Pi(S) \;=\; \mathrm{span}\,(U_S) \;\in\; \Gr(\rho, D),
  \label{eq:grasspoint}
\end{equation}
and the singular values are not retained as coordinates. The truncation
is unique when the $\rho$th singular value is separated from the next one. This is the representation of
\citet{WangWei2025}. Once the normalisation constants are fixed, this map assigns each sequence
a point in a common Grassmann manifold. The distance is then computed between
subspaces rather than between feature vectors.

\begin{definition}[Chordal distance]
Let $\theta_1,\dots,\theta_\rho$ be the principal angles between
$\Pi(S)$ and $\Pi(R)$ \citep[Section 6.4.3]{GolubVanLoan2013}.
Following \citet{Conway1996}, their chordal distance is
\begin{equation}
  \dist_{\mathrm B}(S,R)=\dGr(S,R)
  \;=\;\left(\sum_{i=1}^{\rho}\sin^{2}\theta_i\right)^{1/2}.
  \label{eq:chordal}
\end{equation}
\end{definition}
Since the singular values of $U_S^\top U_R$ are $\cos\theta_i$, an
equivalent expression for computation is
\[
  \dGr(S,R)^2
  =\rho-\bigl\lVert U_S^\top U_R\bigr\rVert_F^2.
\]

\begin{samepage}
\begin{proposition}[The chordal distance is a metric on $\Gr(\rho,D)$]
\label{prop:chordal}
Let $P_S = U_S U_S^{\top}$ denote the orthogonal projector onto $\Pi(S)$, and
let $\theta_1, \dots, \theta_{\rho} \in [0,\pi/2]$ be the principal angles
between $\Pi(S)$ and $\Pi(R)$. Then
\begin{enumerate}
  \item[(i)] $\dGr$ is invariant under $U_S \mapsto U_S O_1$ and
  $U_R \mapsto U_R O_2$ for any $O_1, O_2 \in O(\rho)$, so it is a
  well-defined function on $\Gr(\rho,D) \times \Gr(\rho,D)$;
  \item[(ii)] $\dGr(S,R) = \tfrac{1}{\sqrt{2}}
  \lVert P_S - P_R \rVert_{F}$, and $\dGr$ is a metric on $\Gr(\rho,D)$;
  \item[(iii)]
  $0 \leq \dGr \leq \sqrt{\rho}$, with $\dGr = 0$ iff
  $\Pi(S) = \Pi(R)$ and $\dGr = \sqrt{\rho}$ iff the subspaces are orthogonal.
\end{enumerate}
\end{proposition}

These standard projector and principal-angle expressions for the chordal metric are described by \citet{Conway1996}; see also \citet{Edelman1998}.
\end{samepage}

The projector expression turns subspace comparison into a Euclidean distance
between matrices; the angle expression describes the same distance through
relative orientation. For rank-one subspaces, it reduces to $\dGr=\sin\theta$. Neither
expression depends on which orthonormal basis represents either subspace.

\label{sec:whatdiffers}
The two constructions compare different aspects of the sequence representation.
The Betti-curve distance retains differences in feature amplitudes, whereas
the Grassmann distance compares the orientations of leading subspaces and
also incorporates non-harmonic spectral channels. For a well-defined
truncation, uniform scaling $Z_S\mapsto cZ_S$ with $c>0$ leaves the Grassmann
representation unchanged; reweighting individual channels can change it.
Both constructions admit different alphabets, but their comparison also
depends on the channels, radius grids and normalisation conventions.
The chordal distance is a metric on the Grassmann manifold; on sequences,
distinct inputs can share a representation and hence have zero distance.
Thus the distinction is between comparing measured values in a feature
space and comparing their leading subspaces; the present implementations
also differ in the spectral measurements they include.

\section{Data and methods}

\begin{figure}[!htb]
  \centering
  \includegraphics[width=\textwidth,height=0.70\textheight,keepaspectratio]{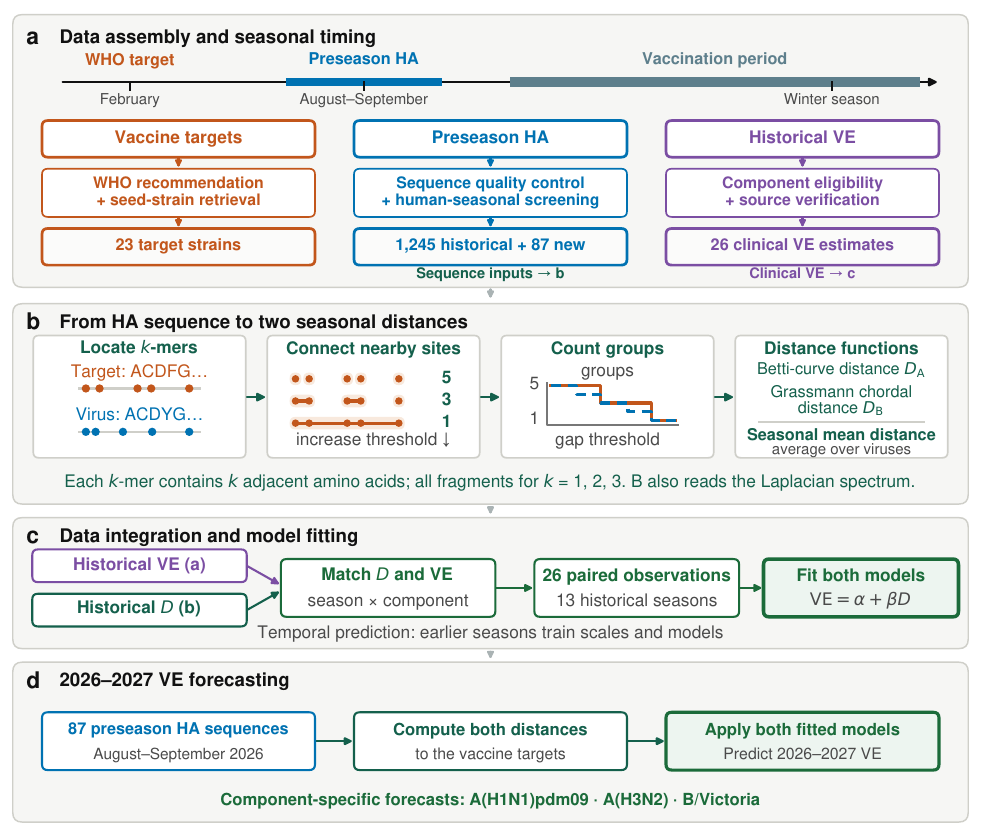}
  \caption{\textbf{Data and prediction procedure.}
  \textbf{(a)} Vaccine targets, preseason HA sequences and published
  component-specific VE estimates are assembled and screened.
  \textbf{(b)} $k$-mer position filtrations provide Betti-0 curves and Laplacian
  descriptors for the two distances, which are averaged by season and component.
  Sequence examples are schematic.
  \textbf{(c)} Distances are paired with 26 historical VE observations to fit
  two regressions; temporal validation fits scales and models on earlier
  seasons only.
  \textbf{(d)} The fitted models use distances from 87 August--September 2026 HA
  sequences to forecast 2026--2027 component-specific VE.}
  \label{fig:dataflow}
\end{figure}

We combined HA sequences collected in the United States during August and
September with historical clinical VE estimates to predict vaccine
effectiveness for the upcoming influenza season. After
sequence filtering, exclusion of nonseasonal viruses and verification of the
published VE estimates, the final historical analysis included 1{,}226 HA
sequences contributing to 26 observations (Supplementary Table~\ref{tab:cohort})
from 13 seasons between 2010--2011
and 2024--2025. Each observation paired the mean distance to the corresponding
vaccine strain with one published VE estimate for a particular season and
viral subtype or lineage. These 26 observations were used to fit the pooled
regressions and to test predictions using earlier seasons to predict later
ones. For the 2026--2027 forecasts, a separate set of 87 HA sequences collected
in August--September 2026 passed the sequence filters: 44 A(H1N1)pdm09,
26 A(H3N2) and 17 B/Victoria. Their mean distances were used to generate
forecasts from the historical regressions; no 2026--2027 VE outcome was used
to fit or test the models. Figure~\ref{fig:dataflow} shows the procedure.

Each observation in the regression contains one mean distance and one VE
estimate for a particular season and viral subtype or lineage. Individual
sequences contribute to the mean distance but do not provide separate clinical
outcomes. This distinction prevents a season with many submitted sequences
from being counted as many independent VE observations.

\subsection{HA sequences and vaccine strains}

We obtained historical HA sequences from GISAID
\citep{ShuMcCauley2017,Elbe2017,Khare2021}. When downloading the data, we set the GISAID collection-date filter to
August 1 through September 30 of each year, inclusive. The searches also required a United States location,
a human host label and an available HA segment. This date filter refers to
when the specimens were collected, rather than when the sequences were
submitted or downloaded. The final historical analysis covered 13 seasons between
2010--2011 and 2024--2025. Supplementary Table~\ref{tab:episets} lists the
GISAID dataset identifiers and DOIs for the searches.

We retained original-passage records with complete HA sequences of at least
1,650 nucleotides for influenza A and at least 1,700 nucleotides for influenza B, and
excluded sequences containing ambiguous bases. We also checked for viruses
that did not represent seasonal human influenza, because the host label
entered by a submitting laboratory was not sufficient to establish this.
Three rules, set before examining VE values, excluded records with a clade
beginning with \texttt{1990}, records assigned to clade \texttt{2010.1}, or
records with less than 95\% gapped-alignment similarity to the corresponding
vaccine strain. The first two rules identify swine-origin lineages; the
sequence comparison also checks records with missing or unassigned clade
labels. Each excluded accession and its reason were retained. These checks
left 1{,}245 historical sequences, of which 1{,}226 belonged to the final
26 observations used in the VE analysis.

For comparison with the circulating viruses, we selected the vaccine strain
recommended for the Northern Hemisphere in each season and for each viral
subtype or lineage \citep{WHOrecommendations}. We obtained the corresponding
HA sequences and checked strain names, sequence sources and passage histories
before calculating distances. Supplementary Table~\ref{tab:targets} lists the
vaccine strains and their sequence accessions.

We used the same search criteria and sequence filters to obtain the data for
the 2026--2027 forecasts. These specimens were collected in the United States
during August--September 2026. The searches were performed on 9 October 2026,
after the collection period ended, and retained 87 HA sequences:
44 A(H1N1)pdm09, 26 A(H3N2) and 17 B/Victoria. We compared these sequences
with the corresponding 2026--2027 vaccine strains and used their mean
distances as inputs to the historical regression models. The GISAID dataset
identifiers and DOIs for these searches are also listed in Supplementary
Table~\ref{tab:episets}.

\subsection{Calculation of sequence distances}
\label{sec:instantiation}

For the distance calculations, we used the HA protein of each circulating
virus and vaccine strain, represented by its ordered sequence of amino acids.
Both methods considered the 20 standard amino acids and recorded the positions
of individual amino acids, adjacent pairs and adjacent triplets
($k=1,2,3$). We calculated the two distances defined in
Section~\ref{sec:framework} between each retained virus and the vaccine strain
for the same season and subtype or lineage. For each method, we then averaged
these virus--vaccine distances within each season and subtype or lineage.
This average was the distance used to predict VE.
Table~\ref{tab:instantiation} lists the settings for both methods.

\begin{table}[htbp]
  \centering
  \caption{Settings for the two sequence distances. Both use the same HA
  sequences and vaccine strains. The listed scales apply to the historical
  pooled fits and frozen forecasts; prediction tests on later seasons estimate
  scales from earlier seasons only.}
  \label{tab:instantiation}
  \small
  \begin{tabular}{lll}
    \toprule
    & \textbf{Topological distance} $\dist_{\mathrm{A}}$
    & \textbf{Geometric distance} $\dist_{\mathrm{B}}$ \\
    & ($\dist_{\mathrm A}$, \citealp{HozumiWei2024})
    & ($\dist_{\mathrm B}$, \citealp{WangWei2025}) \\
    \midrule
    Alphabet            & 20 amino acids      & 20 amino acids \\
    Orders $k$          & $1,2,3$             & $1,2,3$ \\
    Channels $\chan$    & $\beta_0$           & $\beta_0,\ \lambda^{+}_{\min},\ \mu,\ \sigma$ \\
    Radii               & $\{1,2,4,\dots,N\}$, log$_2$-spaced
                        & $T = 50$, log-spaced on $[1,600]$ \\
    Order weights       & $2^{-(K-k)}$ in \eqref{eq:kmer}
                        & $w_k = 2^{-(K-k)}$ in \eqref{eq:Z} \\
    Scale               & global $\omega_k$ (one per order)
                        & global RMS $s_{g,k}$ (one per channel-order) \\
    Rank $\rho$         & ---                 & $10$ \\
    Comparison          & Euclidean, pairwise & chordal on $\Gr(10,D)$ \\
    \bottomrule
  \end{tabular}
\end{table}

The Betti-curve distance counted the connected groups of occurrences of each
amino acid, pair or triplet at radii
$r\in\{1,2,4,8,\dots,N\}$, where $N$ is the greater of the two sequence
lengths. The Grassmann distance used these group counts together with the three
Laplacian summaries defined in Section~\ref{sec:framework}: the smallest
positive eigenvalue, the mean of the positive eigenvalues and their standard
deviation. It used 50 logarithmically spaced radii on $[1,600]$ and retained
10 leading singular vectors to represent each HA protein ($\rho=10$).
The upper radius exceeds the longest HA protein in the data, which contains
585 amino acids. \citet{HozumiWei2024} sample radii linearly as
$r_i=4^{\,k-1}i$, the factor $4^{\,k-1}$ keeping the number of filtration steps
in proportion to the $4^k$ nucleotide $k$-mer types. With the 20-letter
amino-acid alphabet the corresponding factor would advance the $k=3$ radii in
steps of several hundred residues, exceeding the length of an HA protein. We
therefore sample logarithmically, which covers both short and long gaps along
these proteins in a comparable number of steps. We considered individual amino acids, pairs and triplets because the
$20^4=160{,}000$ possible combinations of four amino acids greatly outnumber
the positions in an HA protein.
The protein examples in \citet{HozumiWei2024} also use $k=3$.

We rescaled the measurements before combining the different word lengths
and Laplacian summaries, so that measurements with larger numerical values
would not dominate the distance calculations.
For the historical pooled fits and the frozen forecasts, the Betti-curve
scales were $\omega_1=7.969$, $\omega_2=17.479$ and $\omega_3=24.826$,
calculated across the study comparisons. The Grassmann RMS scales $s_{g,k}$
were fitted using all 529 distinct protein sequences, including those from
2026. Each resulting scale was shared across seasons and viral groups.
Estimating a separate $\omega_k$ for each group would force its mean normalised
distance to a fixed value and remove the differences needed for regression.
The shared scales preserve those differences. Tests of prediction on later
seasons used scales fitted only to earlier data, as described in the next
subsection; the project-wide scales were not used for those tests.

We checked the Betti-curve calculation against the published formulas and
internal consistency tests. For the Grassmann calculation, we used the
reference implementation of \citet{WangWei2025} for the features,
normalisation, SVD and distance, using the 20 amino acids, word lengths of one to three and the
50 radii specified above. With the original nucleotide settings, the calculated feature blocks
matched the reference results to machine precision. Further details of the
settings and implementation checks are given in Supplementary Information,
Section~\ref{sec:implementation}.

\subsection{VE prediction and evaluation}
\label{sec:veaudit}

We matched the mean virus--vaccine distances for each season and subtype or
lineage with the corresponding published clinical VE estimates
\citep{Treanor2012,Ohmit2014,McLean2015,Zimmerman2016,Jackson2017,%
Flannery2019,Rolfes2019,Chung2020,Tenforde2021,Price2023,Chung2024,Chung2025,%
Chung2026}. These estimates came primarily from the US Flu VE Network and
used a test-negative design, which compares vaccination status between
patients who test positive and those who test negative for influenza
\citep{JacksonNelson2013}. All of these estimates are adjusted estimates,
computed as $(1-\mathrm{aOR})\times 100$ from multivariable logistic regression,
where $\mathrm{aOR}$ is the adjusted odds ratio of testing positive for
influenza among vaccinated compared with unvaccinated participants. The source
articles adjust a priori for study site, participant age, presence of
underlying health conditions or general health status, and month of illness
onset, with some seasons additionally adjusting for race and ethnicity, so the
covariate sets are similar but not identical across seasons. The quantity
modelled here, and therefore the quantity forecast, is adjusted
component-specific VE rather than any crude or unadjusted effectiveness.
We checked each estimate against its source article
and required it to match the subtype or lineage being analysed. The final
26 observations, including both mean distances, VE estimates, reported
confidence intervals and source articles, are listed in Supporting
Information, Section~\ref{sec:modeldata} and Supplementary Table~\ref{tab:cohort}. The individual source checks, reasons for excluding
two observations and checks of a figure-derived value are described in
Supplementary Information, Section~\ref{sec:vechecks}.

For each distance, we fitted a separate pooled linear regression,
\begin{equation}
  \VE_i = \alpha + \beta \dist_i + \epsilon_i ,
  \label{eq:pooled}
\end{equation}
where $\dist_i$ is the mean virus--vaccine distance for observation $i$,
$\VE_i$ is the corresponding clinical VE in percentage points, $\alpha$ is
the intercept, $\beta$ is the slope and $\epsilon_i$ is the residual error.
All subtypes and lineages entered the same regression with a common
intercept and slope. We used ordinary least squares, giving each observation
equal weight; sequence counts do not measure the precision of clinical VE.

We assessed the distance--VE association using regression standard errors
and permutations of the VE labels. We report 95\% OLS prediction intervals
for future clinical VE under the regression model. Supplementary Information,
Section~\ref{sec:statdetails}, describes the slope-permutation tests and
prediction intervals.

We tested prediction by fitting the distance scales and regressions to earlier
seasons and predicting the next test season. No sequence from the test season
or a later season was used to estimate those scales. Both models were tested
on the same observations and compared with a baseline that predicted the
mean VE for the relevant subtype or lineage in earlier seasons. We measured
absolute and squared prediction errors and the coverage and width of
95\% prediction intervals. This retrospective evaluation follows season
order rather than reconstructing historical
publication and sequence-deposition dates. Supplementary Information,
Section~\ref{sec:temporaldetails}, gives the test-season eligibility rules,
distance-scale calculations and error measures.

The frozen 2026--2027 forecasts use the original project-wide scales and
historical pooled regressions; they were not replaced by the later evaluation
procedure. We will compare them with final US Flu VE Network estimates for
the matching subtypes and lineage. As a separate check of the preseason
measurements, we compared topological distances from August--September with
those from the subsequent clinical enrollment periods for 18 matched
observations across eight seasons. This comparison did not use VE values;
its methods and data limitations are given in Supplementary Information,
Section~\ref{sec:windowdetails}.

\section{Results}

For the 2026--2027 influenza season in the United States, the pooled
topological model predicted VE of 40.3\% against A(H1N1)pdm09, 56.9\% against
A(H3N2) and 47.1\% against B/Victoria. The pooled geometric model predicted
34.0\%, 56.8\% and 51.1\%, respectively (Table~\ref{tab:forecast}). These
are subtype- and lineage-specific point forecasts, not an overall VE estimate
for the season, and their 95\% prediction intervals are wide. Below, we first
present the historical distance--VE fits used to obtain these forecasts,
then compare August--September topological distances with distances measured
during clinical enrollment, and report the new-season forecasts and their
uncertainty. Prediction tests on later historical seasons are reported in
Supplementary Information, Section~\ref{sec:pooledtemporal}.

\subsection{Pooled relationships between sequence distance and VE}
\label{sec:topological}
\label{sec:geometric}

\begin{figure}[!htb]
  \centering
  \captionsetup{skip=4pt}
  \includegraphics[width=\textwidth]{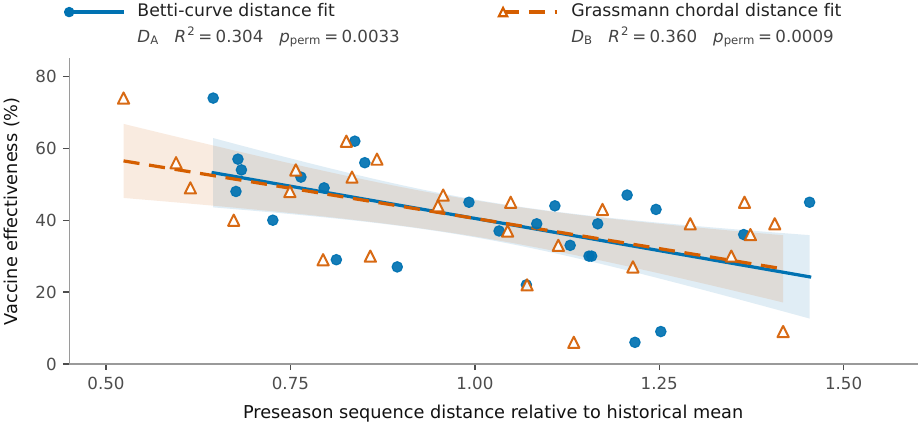}
  \caption{\textbf{Pooled distance--VE fits.} The 26 historical observations
  are shown for both methods: blue circles and a solid line for the Betti-curve
  distance $D_A$; orange open triangles and a dashed line for the Grassmann
  chordal distance $D_B$. Shading shows 95\% confidence bands for mean VE.
  For display, each distance is divided by its historical mean (1 on the
  horizontal axis); this leaves fitted VE, $R^2$ and $p$ values unchanged.
  The slope-permutation $p$ values test each association, not a difference
  between methods. These fits do not measure future prediction accuracy.}
  \label{fig:pooled}
\end{figure}

The historical analysis included 26 observations from 13 seasons, with both
mean virus--vaccine distances paired with the same clinical VE estimate for
each season and subtype or lineage (Supplementary Table~\ref{tab:cohort}).
Larger distances were associated with lower VE for both methods
(Figure~\ref{fig:pooled}). The fitted topological model was
\begin{equation}
  \widehat{\VE}_{\mathrm A}=76.35-22.84\,\dist_{\mathrm A},
\end{equation}
with $R^2=0.304$ and a two-sided slope-permutation $p=0.0033$.
The fitted geometric model was
\begin{equation}
  \widehat{\VE}_{\mathrm B}=74.10-31.53\,\dist_{\mathrm B},
\end{equation}
with $R^2=0.360$ and a two-sided slope-permutation $p=0.0009$.
Thus, the models explained 30.4\% and 36.0\% of the variation in historical
VE, respectively. These pooled fits describe associations across subtypes,
lineages and seasons; they do not establish prediction accuracy for a
future season.

The effects of excluding sequence records and VE observations on the
pooled fit are reported in Supplementary Information, Section~\ref{sec:cohortresults}.
Component-adjusted analyses did not establish
an association beyond subtype or lineage identity (Supplementary Information,
Section~\ref{sec:adjustedspec}). Throughout, \emph{component-adjusted} refers
to the inclusion of subtype or lineage intercepts in the regression, which is
distinct from the confounder adjustment already applied within each clinical VE
estimate (Section~\ref{sec:veaudit}). Separate fits within each subtype or lineage
were also inconclusive (Section~\ref{sec:withingroup}; Figure~\ref{fig:subtypes}).

\subsection{August--September topological distances track enrollment-period distances}
\label{sec:windowresults}

To use August--September sequences for an early forecast, we first examined
whether their distances to the vaccine strain reflected those measured
during the later clinical study. We compared topological distances for 18
matched observations from eight influenza seasons: 2014--2015, 2015--2016,
2016--2017, 2017--2018, 2018--2019, 2019--2020, 2021--2022 and 2022--2023.
Each observation paired the mean distance in August--September with that
during enrollment for one season and subtype or lineage. The comparison
included five A(H1N1)pdm09, six A(H3N2), four B/Victoria and three B/Yamagata
observations. Seasons contributed one to four observations, depending on
which viral groups had sequences available in both windows. Each pair used
the same vaccine strain and shared distance scales. This allowed us to compare the
earlier and later virus samples without changing the vaccine reference
or the scale of the measurements.

\begin{figure}[h!]
  \centering
  \includegraphics[width=\textwidth]{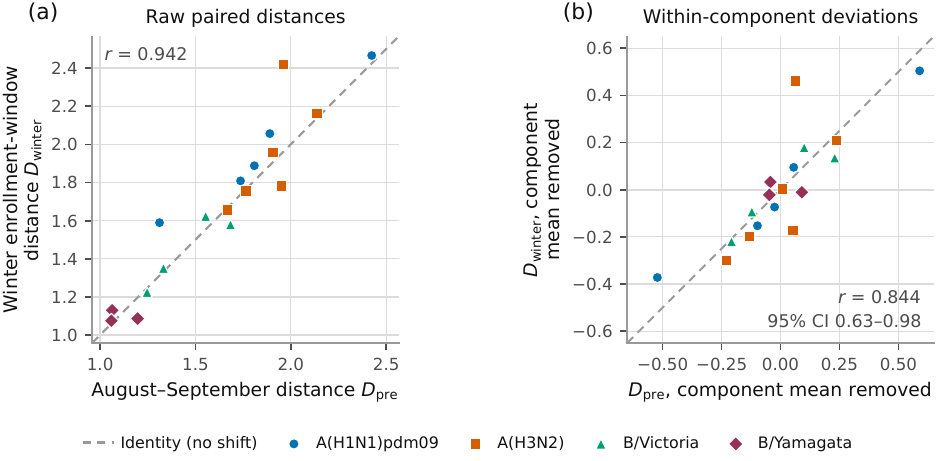}
  \caption{\textbf{Preseason distance tracks winter enrollment-window
  distance.} Both panels show the same 18 season--component observations
  from eight influenza seasons listed in Section~\ref{sec:windowresults}:
  five A(H1N1)pdm09, six A(H3N2), four B/Victoria and three B/Yamagata.
  (a) Each point pairs the mean August--September distance with the mean
  enrollment-period distance for the same season and viral group.
  (b) The same 18 points after subtracting each viral group's mean separately
  in each window, removing between-group differences in mean distance.
  No observations are removed or combined. The dashed line is identity in
  both panels. The 95\% interval in (b) is a season-cluster bootstrap and
  remains a small-sample estimate. This figure establishes temporal concordance
  of the measurement, not its antigenic validity.}
  \label{fig:paired}
\end{figure}

The August--September and enrollment-period distances were strongly
correlated ($r=0.942$; Figure~\ref{fig:paired}a). Observations with larger
virus--vaccine distances before the season generally also had larger
distances during clinical enrollment. However, this pooled correlation
could partly reflect differences between subtypes or lineages that persist
in both windows. We therefore subtracted each viral group's mean separately
in each window and compared the remaining deviations. Both panels show
all 18 observations; subtracting the group means changes their coordinates
without removing or combining observations. The correlation remained high ($r=0.844$; Figure~\ref{fig:paired}b): observations above or
below their group's mean in August--September tended to remain above or below
it during enrollment. The season-cluster bootstrap 95\% interval was
0.632--0.980, reflecting uncertainty across the eight available seasons.

These results support using August--September topological distances as an
early measure of the relative virus--vaccine differences observed during
clinical enrollment. The relationship persisted after removing mean
differences between viral groups, supporting the use of the earlier samples
when enrollment-period sequences are not yet available. Correlation alone
does not establish that the two windows give numerically interchangeable
distances, however. The comparison was limited to seasons with data in both
windows and was performed only for the topological distance. It does not
directly validate the geometric distance or establish accurate clinical VE
prediction. The matching procedure, shared scales and data limitations are
described in Supplementary Information, Section~\ref{sec:windowdetails}.

\subsection{Pooled forecasts for the 2026--2027 season}

The forecast data comprised 44 A(H1N1)pdm09, 26 A(H3N2) and 17 B/Victoria HA
sequences from specimens collected during August--September 2026. The pooled
topological model predicted VE of 40.3\%, 56.9\% and 47.1\%, respectively;
the pooled geometric model predicted 34.0\%, 56.8\% and 51.1\%
(Table~\ref{tab:forecast}). The two methods gave similar point forecasts,
particularly for A(H3N2), and are reported separately without averaging.

The 95\% prediction intervals were approximately 52--58 percentage points wide,
reflecting substantial unexplained variation in the historical fits.
Both sets of forecasts remain frozen for comparison with final 2026--2027
US Flu VE Network estimates for the matching subtypes and lineage.

\begin{table}[hp!]
\centering
\caption{Frozen 2026--2027 component forecasts from the pooled fit of each
distance. The models use two distance representations of the same HA
sequence collection, and their forecasts are reported separately. Intervals are 95\% future-observation intervals.
Both pooled forecasts are frozen and will be evaluated prospectively. No
2026--2027 VE outcome was used in fitting.}
\label{tab:forecast}
\small
\setlength{\tabcolsep}{4pt}
\begin{tabular}{cccccccc}
\toprule
 & & \multicolumn{3}{c}{$k$-mer topological ($D_\mathrm{A}$)} & \multicolumn{3}{c}{Grassmann geometric ($D_\mathrm{B}$)} \\
\cmidrule(lr){3-5}\cmidrule(lr){6-8}
Component & $n$ viruses & Mean $D$ & \shortstack{2026--2027\\Predicted VE (\%)} & 95\% PI & Mean $D$ & \shortstack{2026--2027\\Predicted VE (\%)} & 95\% PI \\
\midrule
A(H1N1)pdm09 & 44 & 1.5804 & 40.3 & [13, 68] & 1.2706 & 34.0 & [8, 60] \\
A(H3N2) & 26 & 0.8506 & 56.9 & [28, 86] & 0.5497 & 56.8 & [29, 84] \\
B/Victoria & 17 & 1.2814 & 47.1 & [20, 75] & 0.7301 & 51.1 & [24, 78] \\
\bottomrule
\end{tabular}

\end{table}

\section{Discussion}

This study linked HA sequence distances to historical clinical VE and used
the fitted relationships to forecast subtype- and lineage-specific VE for
the 2026--2027 US influenza season. Both distances showed pooled associations
with VE, and August--September topological distances tracked those measured
during later clinical enrollment. Together, these findings support
investigating preseason sequences as an early source of information for
VE prediction. The limited improvement in historical prediction and wide
forecast intervals nevertheless leave practical reliability to be
established through prospective evaluation.

\medskip
\noindent\textbf{Comparability of overall and component-specific VE estimates.}
The CDC past-seasons summary reports a single overall VE estimate against
any influenza across all ages~\cite{CDCPastSeasons}. The forecasts here,
however, concern individual subtypes and lineages.
Both are adjusted clinical estimates in the sense of
Section~\ref{sec:veaudit}; the distinction is whether the outcome combines
co-circulating viruses or refers to a particular viral group.
This distinction matters when interpreting the forecasts and comparing
them with subsequently reported clinical VE.

Table~\ref{tab:cdccomparison} compares the overall and component-specific
estimates across the 13 seasons represented in this study.
CDC overall VE exceeded the unweighted mean of the available component
estimates by an average of 2.4 percentage points.
Excluding 2024--2025 reduced this average difference to 0.8 percentage
points. These small average differences conceal variation among viral
groups. Overall VE exceeded the A(H3N2) estimate in 10 of 11 seasons
and fell below the available B estimates in seven of eight seasons.
This pattern is consistent with the lower VE commonly reported against
A(H3N2)~\cite{Belongia2016}, but the available components do not represent
every circulating group in every season. In particular, the four earliest
seasons contributed only A(H3N2) observations.
Their component means therefore equal the A(H3N2) estimates and cannot
represent overall seasonal VE. The two largest negative differences, $-7$ in 2016--2017 and $-17$ in
2022--2023, coincide with A(H3N2) predominance and reduced vaccine match from
egg-adapted vaccine strains in the former~\cite{Flannery2019,Zost2017}, and
with an early season following minimal influenza circulation during the
COVID-19 pandemic in the latter~\cite{Olsen2021,Chung2024}.

Differences in surveillance network also affect the comparison.
The 2024--2025 overall entry on the CDC summary page comes from the
VISION network, whereas the component estimates used here come from
the US Flu VE Network~\cite{CDC20242025,Frutos2025}.
The overall estimates from the two networks differed by 14 percentage
points in that season. Such differences reinforce the need to match
the clinical setting, study population and viral group when evaluating
a forecast.

\begin{table}[!htbp]
\centering
\caption{The CDC overall VE estimate for each season compared with the
component-specific estimates used in this study.
CDC values are adjusted overall VE (CDC OV) against any influenza
across all ages, with 95\% confidence intervals, as listed on the
CDC past-seasons page~\cite{CDCPastSeasons}.
Component values are the clinical estimates in Table~\ref{tab:cohort};
a dash indicates that no observation for that group was included
for the season.
The mean of the components (MC) is the unweighted average of the
available group estimates and should not be interpreted as an
estimate of overall seasonal VE.
Differences are expressed in percentage points.
$^{a}$The 2024--2025 CDC entry comes from the VISION network,
whereas the component estimates come from the US Flu VE
Network~\cite{CDC20242025,Frutos2025}.}
\label{tab:cdccomparison}
\small
\setlength{\tabcolsep}{4pt}
\begin{tabular}{lccccccc}
\toprule
 & \shortstack{CDC adjusted\\overall VE (\%)} & \multicolumn{4}{c}{Component-specific VE (\%)} & & Difference \\
\cmidrule(lr){2-2}\cmidrule(lr){3-6}\cmidrule(lr){8-8}
Season & \shortstack{All ages, any\\influenza (CDC OV)} & H1N1 & H3N2 & B/Vic & B/Yam & \shortstack{Mean of\\components (MC)} & \shortstack{CDC OV\\$-$ MC} \\
\midrule
2010--2011 & 60 [53, 66] & \textemdash & 54 & \textemdash & \textemdash & 54.0 & +6.0 \\
2011--2012 & 47 [36, 56] & \textemdash & 39 & \textemdash & \textemdash & 39.0 & +8.0 \\
2012--2013 & 49 [43, 55] & \textemdash & 39 & \textemdash & \textemdash & 39.0 & +10.0 \\
2014--2015 & 19 [10, 27] & \textemdash & 6 & \textemdash & \textemdash & 6.0 & +13.0 \\
2015--2016 & 48 [41, 55] & 45 & 43 & 49 & 57 & 48.5 & $-$0.5 \\
2016--2017 & 40 [32, 46] & \textemdash & 33 & 56 & 52 & 47.0 & $-$7.0 \\
2017--2018 & 38 [31, 43] & 62 & 22 & \textemdash & 48 & 44.0 & $-$6.0 \\
2018--2019 & 29 [21, 35] & 44 & 9 & \textemdash & \textemdash & 26.5 & +2.5 \\
2019--2020 & 39 [32, 44] & 30 & \textemdash & 45 & \textemdash & 37.5 & +1.5 \\
2021--2022 & 36 [21, 48] & \textemdash & 36 & \textemdash & \textemdash & 36.0 & +0.0 \\
2022--2023 & 30 [$-$9, 54] & 47 & \textemdash & \textemdash & \textemdash & 47.0 & $-$17.0 \\
2023--2024 & 44 [36, 51] & 29 & 30 & 74 & \textemdash & 44.3 & $-$0.3 \\
2024--2025 & 56 [54, 58]$^{a}$ & 37 & 27 & 40 & \textemdash & 34.7 & +21.3 \\
\midrule
\multicolumn{2}{l}{Mean difference, CDC OV $-$ component} & $-$1.4 & +11.6 & $-$7.4 & $-$10.3 & & +2.4 \\
\multicolumn{2}{l}{Seasons with CDC OV higher} & 4/7 & 10/11 & 1/5 & 0/3 & & 7/13 \\
\bottomrule
\end{tabular}

\end{table}

We therefore specify the intended clinical comparison in advance.
The frozen forecasts are to be evaluated against final US Flu VE Network
adjusted estimates for medically attended outpatient illness with
laboratory-confirmed influenza, separately for A(H1N1)pdm09, A(H3N2)
and B/Victoria. Evaluation should consider both the point predictions
and their prediction intervals. Comparing these forecasts directly with
a single overall VE estimate would not assess the outcomes the models
were fitted to predict.

\medskip
\noindent\textbf{Interpreting the pooled associations.}
The fitted relationships combine differences across subtypes, lineages
and seasons. They do not establish that changes in distance within
each viral group correspond to changes in VE.
After accounting for differences in mean VE among groups, the shared
distance slope was not statistically significant
(Supplementary Information, Section~\ref{sec:adjustedspec};
Supplementary Table~\ref{tab:models}).
Separate fits for individual subtypes and lineages were also inconclusive
(Section~\ref{sec:withingroup}; Figure~\ref{fig:subtypes}).
The pooled association may therefore reflect both between-group
differences and seasonal variation within groups, whose contributions
remain unresolved. The small number of observations within each group
also limits the interpretation of nonsignificant slopes.

The window comparison addresses a different question: whether earlier
samples reflect relative virus--vaccine differences during clinical
enrollment. Its positive result supports the timing of the measurement
in the available seasons, but does not independently validate the
relationship with clinical VE.
Agreement varied across groups and seasons and was assessed only for
the topological distance. Detailed influence checks and group-specific
results are provided in Supplementary Information,
Section~\ref{sec:discussiondetails}.

\medskip
\noindent\textbf{Historical prediction accuracy and uncertainty.}
Fitting distance scales and regression coefficients using earlier seasons
provides a more direct assessment of forecasting performance than the
pooled fit alone. Across 19 observations from nine test seasons, mean
absolute errors were 12.01 percentage points for the topological model
and 11.34 for the geometric model.
The reference that predicted each component's earlier mean VE had an
error of 11.96 percentage points.
The geometric model's improvement over this reference was therefore
0.62 percentage points. Mean absolute error and root mean squared error
also ranked the two distance models differently.
These results indicate limited improvement over the reference and do not
demonstrate a consistent advantage for either distance representation.

Prediction intervals further indicate substantial uncertainty around
individual VE forecasts. Their coverage must be considered together with
their width, because wide intervals can contain observed VE while offering
limited precision.
Prediction errors, interval coverage and interval widths are reported in
Supplementary Information, Section~\ref{sec:pooledtemporal} and
Supplementary Table~\ref{tab:validationall}.
The evaluation followed season order without fully reconstructing the
information available at historical forecast dates.
Observations from the same season also share a fitted model and should
not be treated as independent seasonal replications.

The two methods differ in their sequence summaries, radius settings and
normalisation, as described in Section~\ref{sec:implementation}.
Their performance differences therefore cannot be attributed to
Grassmann geometry alone.
Prospective evaluation against the specified clinical outcomes will
help determine whether either representation provides useful predictive
information beyond historical VE.

\section{Conclusion}

Topological and geometric distances between August--September HA sequences
and vaccine strains were associated with historical clinical VE across
26 season--component observations. The topological distances also tracked
those measured during later clinical enrollment in the available paired
seasons, supporting their use as an early measure of relative virus--vaccine
differences. Using the historical pooled relationships, we generated frozen
2026--2027 US VE forecasts of 40.3\%, 56.9\% and 47.1\% for A(H1N1)pdm09,
A(H3N2) and B/Victoria, respectively, with the topological model, and
34.0\%, 56.8\% and 51.1\% with the geometric model. These are component-specific
point forecasts with wide prediction intervals. Historical prediction gains
were limited, and information beyond subtype or lineage identity was not well
established. Comparison with the corresponding final clinical VE estimates,
together with evaluation on additional independent data, is needed to
determine whether these sequence-based forecasts can reliably support
early-season planning.

\section*{Data and code availability}

Analysis code, record-level exclusion reasons with accessions, aggregated model
inputs, and the saved forecasts and figure data are available at
\url{https://github.com/XiangXiangJY/Antigenic}.
GISAID sequences are subject to GISAID data-access terms and are identified by
EPI\_SET accessions and DOIs (Supplementary Table~\ref{tab:episets}) rather
than redistributed.
Every preseason search contributing to the cohort or the forecast carries an
EPI\_SET identifier (Supplementary Table~\ref{tab:episets}). The 2022--2023
A(H1N1)pdm09 identifier was generated on 14 October 2026 by re-running the
original search specification, which returned exactly the same 99 records as the
original download.

\section*{Acknowledgements}

We acknowledge the authors, originating laboratories and submitting
laboratories that contributed the specimens, sequences and metadata shared
through GISAID (Supplementary Table~\ref{tab:episets}).
This work was supported in part by NIH grant R01AI164266, the University of Georgia (UGA), and the Georgia Research Alliance.

\bibliographystyle{unsrtnat}
\bibliography{references}

\appendix
\section*{Supplementary Information}
\renewcommand{\thetable}{S\arabic{table}}
\setcounter{table}{0}
\renewcommand{\thefigure}{S\arabic{figure}}
\setcounter{figure}{0}
\section{Component-adjusted and within-group associations}
\label{sec:component}

The pooled fits combine observations from different subtypes and lineages.
Here, we examine whether distance remains associated with VE after allowing
separate group intercepts, and whether a relationship is evident within each
group separately. Neither analysis established a significant distance--VE association.
These results limit the interpretation of the pooled associations, while
the small group sizes leave the within-group relationships uncertain.

\subsection{Different subtype and lineage intercepts with a shared distance slope}
\label{sec:adjustedspec}
\label{sec:adjustedresults}

To assess whether distance provided information beyond subtype or lineage,
we fitted
\begin{equation}
  \VE_i=\alpha_{s(i)}+\beta\dist_i+\epsilon_i,
  \label{eq:subtype}
\end{equation}
where $s(i)$ identifies the subtype or lineage of observation $i$.
Each of A(H1N1)pdm09, A(H3N2), B/Victoria and B/Yamagata has a separate
intercept $\alpha_{s(i)}$, while $\beta$ is a single distance slope shared
by all four groups. Thus, the model allows different VE levels between
groups and estimates the distance--VE relationship after accounting for
these mean differences. The topological slope was
$-13.40$ (SE 9.32; 95\% CI, $-32.8$ to $6.0$; $p=0.1651$; model
$R^2=0.4420$; Table~\ref{tab:models}). The geometric slope was $-17.43$
(SE 12.11; 95\% CI, $-42.6$ to $7.8$; $p=0.1649$; model $R^2=0.4421$). The two
adjusted fits happen to agree closely on both quantities.
Neither adjusted slope was significant.

Adjustment reduced the magnitude of the topological slope from $-22.84$
to $-13.40$ and the geometric slope from $-31.53$ to $-17.43$.
This attenuation is consistent with part of the pooled association reflecting
differences between viral groups. Both adjusted confidence intervals
included zero, so these data did not establish that distance contributed
information about VE beyond subtype or lineage identity. However, the wide
intervals also allow substantial negative slopes; they do not demonstrate
that the within-group relationship is absent.

The adjusted models had higher $R^2$ than the pooled models, but the extra
group intercepts also allow a closer fit to the same observations.
This increase describes overall fit and does not establish an additional
contribution from distance or improved prediction on later seasons.
The common-slope assumption also means that this analysis does not test
whether the distance--VE relationship differs between groups; separate
exploratory fits are described in Section~\ref{sec:withingroup}.

Excluding the figure-derived 2017--2018 B/Yamagata VE changed the adjusted
topological slope from $-13.40$ ($t=-1.44$) to $-13.60$ ($t=-1.43$).
The conclusion was therefore unchanged when this observation was omitted.

\begin{table}[htbp]
\centering
\caption{Topological pooled and component-adjusted regression estimates on
26 observations. Standard errors are ordinary OLS.}
\label{tab:models}
\small
\begin{tabular}{lrrrrr}
\toprule
Term & Estimate & SE & 95\% CI & OLS $p$ & $R^2$ \\
\midrule
\textbf{Primary: pooled} &  &  &  &  &  \\
\quad Intercept & 76.35 & 11.35 &  &  &  \\
\quad Distance $D$ & -22.84 & 7.05 & [-37.4, -8.3] & 0.003 & 0.304 \\
\textbf{Sensitivity: subtype fixed effects} &  &  &  &  &  \\
\quad Distance $D$ & -13.40 & 9.32 & [-32.8, 6.0] & 0.165 & 0.442 \\
\quad Intercept: A(H1N1)pdm09 & 64.86 & 16.57 &  &  &  \\
\quad Intercept: A(H3N2) & 54.18 & 16.73 &  &  &  \\
\quad Intercept: B/Victoria & 69.67 & 12.98 &  &  &  \\
\quad Intercept: B/Yamagata & 67.19 & 12.57 &  &  &  \\
\bottomrule
\end{tabular}

\end{table}

\subsection{Separate distance--VE fits for each subtype and lineage}
\label{sec:withingroup}

We examined whether the pooled distance--VE relationship was also evident
across seasons within individual subtypes or lineages. For this exploratory
analysis, we fitted a separate intercept and slope for each viral group,
rather than imposing the common slope used in Section~\ref{sec:adjustedspec}.
The historical dataset contained seven A(H1N1)pdm09, eleven A(H3N2), five
B/Victoria and three B/Yamagata observations.

\begin{figure}[htbp]
  \centering
  \includegraphics[width=\textwidth]{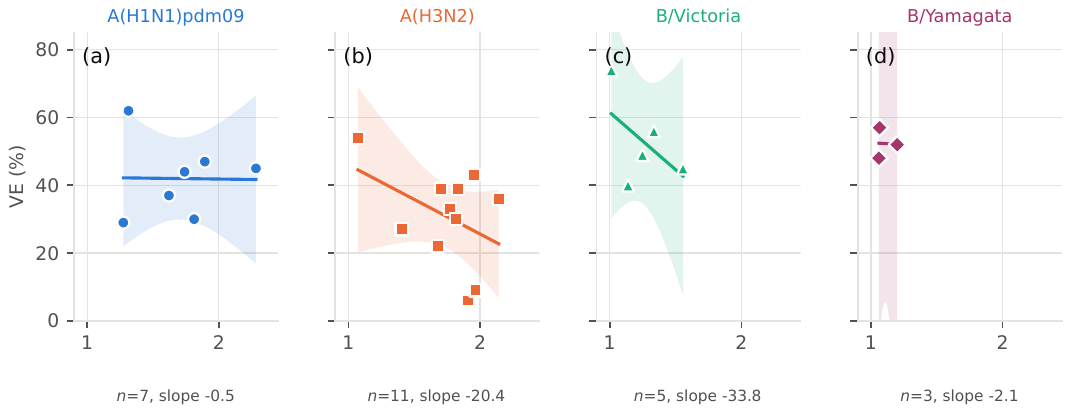}
  \caption{\textbf{Separate distance--VE fits for each subtype and lineage.}
  (a--d) Historical observations for A(H1N1)pdm09 ($n=7$), A(H3N2) ($n=11$),
  B/Victoria ($n=5$) and B/Yamagata ($n=3$), respectively. In every panel,
  the horizontal axis is the mean August--September topological distance
  $D_A$ to the vaccine strain, and the vertical axis is clinical VE (\%).
  Each group has a separately estimated intercept and slope. Solid lines
  and shading show ordinary least-squares fits and 95\% confidence bands
  for mean VE in all four panels. B/Yamagata has only three observations
  spanning a distance range of 0.14, resulting in a very wide confidence band
  that extends beyond the displayed VE range.
  All four slope estimates are negative, but their 95\% confidence intervals
  include zero. These exploratory fits were not used for forecasting.}
  \label{fig:subtypes}
\end{figure}

Within-group analyses did not establish a distance--VE association for any
subtype or lineage with either distance. Figure~\ref{fig:subtypes}a--d shows
the topological results. All four topological slope estimates were negative,
but their 95\% confidence intervals included zero. The small group sizes
and wide intervals leave the within-group relationships unresolved. These
fits assess whether the pooled association also appears within individual
groups and were not used to generate the 2026--2027 forecasts.
In particular, the nearly flat H1N1 estimate and the narrow range of
B/Yamagata distances provide little evidence of a common pattern within
every group. Negative point estimates alone do not establish such a pattern.

\section{Pooled model fits and predictions for later seasons}
\label{sec:additionalresults}

This section examines how excluding sequence records and VE observations
affects the pooled fits, then evaluates predictions for later seasons using
only earlier seasons for training. The fitted relationships changed after
data exclusions. In historical prediction tests, the geometric pooled model
modestly reduced absolute error relative to the historical-mean reference,
but prediction intervals remained wide.

\subsection{How sequence and VE exclusions changed the pooled fit}
\label{sec:cohortresults}
\label{sec:sequenceaudit}

Sequence quality checks left 1,285 historical and 87 forecast-season records.
The seasonal-human-virus checks excluded 40 historical records: 15 by the
clade-prefix rule, 19 by clade 2010.1 and six by the similarity rule when
applied in that order. Independently, the similarity rule flagged all 40.
Of 1,245 retained historical sequences, 1,226 contributed to the final
26 observations (7 H1N1, 11 H3N2, five B/Victoria and three B/Yamagata).
Clinical source checks reduced 28 assembled observations to 26; the two
exclusions were fixed before refitting (Section~\ref{sec:vechecks}). They
changed the pooled topological slope from $-19.26$ ($R^2=0.190$) to
$-22.84$ ($R^2=0.304$).

With the same 26 clinical outcomes, using the sequence sets before the
seasonal-human-virus exclusions gave slope $-9.02$ ($R^2=0.131$;
slope-permutation $p=0.075$), compared with $-22.84$ ($R^2=0.304$;
$p=0.0033$) after exclusions. Six group means changed; the largest changes
were in 2012--2013 H3N2 (11 of 12 records removed) and 2010--2011 H3N2
(two of five removed). These changes show that the pooled fit depends
materially on which sequence records and clinical estimates enter the
analysis. Its interpretation therefore rests on the documented eligibility
rules and source checks. Table~\ref{tab:cohort} provides the final paired data.

\subsection{Historical prediction performance of the pooled models}
\label{sec:pooledtemporal}

We assessed whether the pooled models could predict VE for a later season
using distances and VE observations from earlier seasons. At each test
season, distance scales and regression coefficients were re-estimated from
the earlier seasons only. Both models were evaluated on the same 19
season--component observations across nine test seasons. Each pooled model
was refitted nine times, once per test season. The fit for a season produced
predictions for all eligible viral groups in that season. Each prediction
corresponds to one clinical VE estimate for a particular season and subtype
or lineage. The eligible seasons and inclusion rules are specified in
Section~\ref{sec:temporaldetails}.

Table~\ref{tab:validationall} compares prediction errors with a reference
that uses the earlier mean VE for the relevant subtype or lineage.
This reference asks whether sequence distances add predictive value beyond
previous VE levels. The geometric pooled MAE was 11.34 percentage points and
the topological pooled MAE was 12.01, against 11.96 for the reference. Both
models had lower RMSE than the reference, the topological model the lower of
the two.

These summaries rank the approaches differently. The topological model had the
slightly higher mean absolute error (12.01 against 11.96) yet was closer to the
observed VE in 11 of the 19 tests, against 8 for the reference, a split that an
exact sign test does not distinguish from chance ($p=0.65$); its absolute errors
also had the smallest spread of the three (standard deviation 7.97 percentage
points, against 9.56 for the geometric model and 9.69 for the reference). The
higher mean therefore reflects a few large errors rather than a systematic
difference. None of the paired differences in absolute error was distinguishable
from zero: topological minus reference $+0.05$ percentage points (paired
bootstrap 95\% interval $-2.97$ to $+3.12$), geometric minus reference $-0.62$
($-2.61$ to $+1.52$), and geometric minus topological $-0.67$ ($-3.58$ to
$+2.53$). With 19 tests, a split of at least 15 to 4 would be needed to reach
the conventional significance threshold. These seasons therefore cannot
establish an ordering among the two distances and the reference.

The 95\% prediction intervals contained the observed VE in 19 of 19 tests
for the topological model and 18 of 19 for the geometric model.
Their mean widths were 61.0 and 63.0 percentage points, respectively.
High coverage therefore accompanied wide uncertainty about individual VE
outcomes. Multiple observations came from the same season, and the evidence
spans only nine test seasons. These results describe retrospective prediction
performance and do not establish reliable prospective forecasts.

\begin{table}[htbp]
\centering
\caption{Prediction accuracy and interval coverage for 19 historical VE
observations across nine test seasons. Every row uses the same 19
season--component observations. Each pooled model was refitted for each
test season using earlier seasons only. MAE is mean absolute error;
RMSE is root mean squared error. Both errors and interval widths are in
percentage points (pp). ``Observed VE inside PI'' counts how many actual
VE estimates fall within their predicted 95\% intervals. Coverage is this
count divided by 19, expressed as a percentage. Thus, the geometric model
made 19 predictions, and 18 observed VE estimates fell inside its intervals.}
\label{tab:validationall}
\small
\setlength{\tabcolsep}{3pt}
\begin{tabular}{lcccccc}
\toprule
 & & \multicolumn{2}{c}{Prediction error} & \multicolumn{3}{c}{95\% prediction intervals (PI)} \\
\cmidrule(lr){3-4}\cmidrule(lr){5-7}
Model or reference & \shortstack{Total\\predictions} & \shortstack{MAE\\(pp)} & \shortstack{RMSE\\(pp)} & \shortstack{Observed VE\\inside PI\\(count)} & \shortstack{Coverage\\(\%)} & \shortstack{Mean width\\(pp)} \\
\midrule
Historical-mean reference & 19 & 11.96 & 15.23 & 19 & 100.0 & 77.9 \\
Topological pooled & 19 & 12.01 & 14.30 & 19 & 100.0 & 61.0 \\
Geometric pooled & 19 & 11.34 & 14.67 & 18 & 94.7 & 63.0 \\
\bottomrule
\end{tabular}

\end{table}

Table~\ref{tab:predictiondetail} lists the 19 individual predictions behind
these summaries. Absolute errors ranged from 2.5 to 30.1 percentage points
for the topological model and from 0.9 to 38.3 percentage points for the
geometric model, and both models predicted a higher VE than was observed in
11 of the 19 tests. The geometric prediction was closer to the observed VE
in 13 of the 19 tests, but its single largest error was larger than any
topological error; the two error measures therefore rank the models
differently. Each model was closer to the observed VE than the
historical-mean reference in about half of the tests, 11 of 19 for the
topological model and 10 of 19 for the geometric model, so the small
differences in average error do not correspond to a consistent advantage on
individual season--component outcomes.

The largest error for all three approaches occurred in the 2014--2015
A(H3N2) test, which had the lowest observed VE in the data (6\%) and was
also the earliest test season, fitted on only three earlier observations.
This test is the one case in which an observed VE fell outside a prediction
interval, the lower limit of the geometric interval being 12.3\%. Individual
interval widths ranged from 53.9 to 79.2 percentage points for the
topological model and from 55.4 to 91.4 for the geometric model. Predictions
for single season--component outcomes are therefore compatible with a wide
range of VE values, and the historical-mean reference intervals are wider
still, extending below zero and above 100\% in several early tests.

\begin{table}[htbp]
\centering
\caption{Individual predictions for the 19 historical test observations
summarised in Table~\ref{tab:validationall}. Each row is one clinical VE
estimate for one season and subtype or lineage. ``Training observations''
is the number of earlier season--component observations used to re-estimate
the distance scales and regression coefficients for that test season; all
rows from the same test season share the same fit. Predictions and
intervals are in VE percentage points, and the intervals are 95\%
future-observation intervals from the fit for that season. They are not
constrained to the 0--100\% range, so some reference limits fall outside it.
$^{\dagger}$ marks the one prediction whose interval does not contain the
observed VE. The final rows repeat the summary measures of
Table~\ref{tab:validationall}, recomputed from the rows shown here.}
\label{tab:predictiondetail}
\footnotesize
\setlength{\tabcolsep}{3pt}
\begin{tabular}{llcccccc}
\toprule
 & & & & \multicolumn{3}{c}{Predicted VE (95\% prediction interval)} \\
\cmidrule(lr){5-7}
Test season & \shortstack[l]{Subtype or\\lineage} & \shortstack{Training\\observations} & \shortstack{Observed\\VE (\%)} & \shortstack{Topological\\pooled} & \shortstack{Geometric\\pooled} & \shortstack{Historical-mean\\reference} \\
\midrule
2014--2015 & A(H3N2) & 3 & 6 & 36.1 (3.3, 68.9) & 44.3 (12.3, 76.2)$^{\dagger}$ & 44.0 (1.0, 87.0) \\
2016--2017 & A(H3N2) & 8 & 33 & 38.9 (1.2, 76.6) & 39.7 (0.2, 79.1) & 36.2 ($-$18.4, 90.8) \\
2016--2017 & B/Victoria & 8 & 56 & 47.0 (8.5, 85.5) & 53.5 (7.8, 99.2) & 49.0 ($-$21.5, 119.5) \\
2016--2017 & B/Yamagata & 8 & 52 & 49.5 (9.9, 89.1) & 46.5 (6.0, 86.9) & 57.0 ($-$13.5, 127.5) \\
2017--2018 & A(H1N1)pdm09 & 11 & 62 & 48.6 (19.0, 78.2) & 46.7 (16.5, 77.0) & 45.0 ($-$1.2, 91.2) \\
2017--2018 & A(H3N2) & 11 & 22 & 40.9 (11.8, 70.0) & 40.7 (10.6, 70.8) & 35.7 (0.4, 71.0) \\
2017--2018 & B/Yamagata & 11 & 48 & 53.9 (22.8, 85.0) & 49.3 (18.5, 80.1) & 54.5 (14.5, 94.5) \\
2018--2019 & A(H1N1)pdm09 & 14 & 44 & 38.5 (9.7, 67.3) & 43.0 (13.7, 72.3) & 53.5 (18.2, 88.8) \\
2018--2019 & A(H3N2) & 14 & 9 & 33.4 (3.7, 63.2) & 27.7 ($-$4.9, 60.4) & 33.7 (2.9, 64.5) \\
2019--2020 & A(H1N1)pdm09 & 16 & 30 & 34.7 (4.9, 64.6) & 47.7 (19.1, 76.3) & 50.3 (15.8, 84.9) \\
2019--2020 & B/Victoria & 16 & 45 & 41.5 (12.0, 71.0) & 26.1 ($-$3.8, 56.0) & 52.5 (15.8, 89.2) \\
2021--2022 & A(H3N2) & 18 & 36 & 25.7 ($-$3.4, 54.8) & 28.6 ($-$2.1, 59.3) & 30.6 ($-$0.4, 61.7) \\
2022--2023 & A(H1N1)pdm09 & 19 & 47 & 33.5 (6.3, 60.7) & 43.2 (14.8, 71.7) & 45.2 (13.7, 76.8) \\
2023--2024 & A(H1N1)pdm09 & 20 & 29 & 49.1 (21.7, 76.4) & 47.5 (19.6, 75.4) & 45.6 (15.8, 75.4) \\
2023--2024 & A(H3N2) & 20 & 30 & 36.3 (9.4, 63.2) & 30.9 (2.4, 59.3) & 31.2 (2.5, 59.9) \\
2023--2024 & B/Victoria & 20 & 74 & 55.2 (26.7, 83.7) & 56.5 (26.6, 86.4) & 50.0 (18.6, 81.4) \\
2024--2025 & A(H1N1)pdm09 & 23 & 37 & 40.5 (12.9, 68.2) & 40.2 (12.5, 67.9) & 42.8 (13.0, 72.6) \\
2024--2025 & A(H3N2) & 23 & 27 & 46.0 (18.2, 73.9) & 34.0 (6.0, 62.0) & 31.1 (2.2, 60.0) \\
2024--2025 & B/Victoria & 23 & 40 & 52.8 (24.2, 81.3) & 52.6 (24.1, 81.1) & 56.0 (25.1, 86.9) \\
\midrule
\multicolumn{4}{l}{Mean absolute error (pp)} & 12.01 & 11.34 & 11.96 \\
\multicolumn{4}{l}{Root mean squared error (pp)} & 14.30 & 14.67 & 15.23 \\
\multicolumn{4}{l}{Observed VE inside interval} & 19/19 & 18/19 & 19/19 \\
\multicolumn{4}{l}{Mean interval width (pp)} & 61.0 & 63.0 & 77.9 \\
\bottomrule
\end{tabular}

\end{table}

\section{Sequence records, VE sources and calculation procedures}
\label{sec:methoddetails}

This section documents the data and calculations behind the reported results.
It first lists the 26 distance--VE observations used to fit the models,
followed by vaccine strain records, GISAID sequence sets and checks of the
published VE estimates. The remaining subsections describe distance
calculations, statistical tests, prediction evaluation and the matching of
sampling windows. They also identify limitations in record availability,
clinical estimates and the coverage of the window comparison.

\subsection{Historical distances and VE estimates used to fit the models}
\label{sec:modeldata}

Table~\ref{tab:cohort} lists the 26 observations used to fit the historical
VE models. Each row represents one season and subtype or lineage and pairs
its mean virus--vaccine distance with a published clinical VE estimate.
The two distance columns were calculated from the same sequence records;
the sequence count identifies how many viruses contributed to each mean.
Clinical confidence intervals describe uncertainty in the published VE
estimates but were not used as regression weights. The source checks and
exclusions are explained in Section~\ref{sec:vechecks}; their effect on the
pooled fit is reported in Section~\ref{sec:cohortresults}.

\begin{table}[!htbp]
\centering
\caption{Historical distance--VE data used to fit the models (26 observations). Mean $D_A$ and $D_B$ are
the original-scale Betti-curve and Grassmann distances, rounded to four
decimal places, calculated from the same $n$ viruses against the corresponding
vaccine strain. VE estimates and references identify the clinical outcomes;
95\% CIs are recorded where available and do not enter the regression.
Source-specific qualifications are given in Section~\ref{sec:vechecks}.}
\label{tab:cohort}
\small
\resizebox{\textwidth}{!}{\begin{tabular}{llrrrrcl}
\toprule
Season & Component & $n$ viruses & Mean $D_A$ & Mean $D_B$ & VE (\%) & 95\% CI & VE source \\
\midrule
2010--2011 & A(H3N2) & 3 & 1.0728 & 0.8070 & 54 & [42, 64] & \citep{Treanor2012} \\
2011--2012 & A(H3N2) & 4 & 1.8310 & 1.4990 & 39 & [23, 52] & \citep{Ohmit2014} \\
2012--2013 & A(H3N2) & 1 & 1.7012 & 1.3769 & 39 & [29, 47] & \citep{McLean2015} \\
2014--2015 & A(H3N2) & 10 & 1.9102 & 1.2086 & 6 & [-5, 17] & \citep{Zimmerman2016} \\
2015--2016 & A(H1N1)pdm09 & 8 & 2.2818 & 1.1175 & 45 & --- & \citep{Jackson2017} \\
2015--2016 & A(H3N2) & 209 & 1.9550 & 1.2501 & 43 & --- & \citep{Jackson2017} \\
2015--2016 & B/Victoria & 5 & 1.2487 & 0.6547 & 49 & --- & \citep{Jackson2017} \\
2015--2016 & B/Yamagata & 18 & 1.0654 & 0.9243 & 57 & --- & \citep{Jackson2017} \\
2016--2017 & A(H3N2) & 99 & 1.7723 & 1.1863 & 33 & [23, 41] & \citep{Flannery2019} \\
2016--2017 & B/Victoria & 12 & 1.3353 & 0.6340 & 56 & [23, 75] & \citep{Flannery2019} \\
2016--2017 & B/Yamagata & 12 & 1.1994 & 0.8882 & 52 & [42, 61] & \citep{Flannery2019} \\
2017--2018 & A(H1N1)pdm09 & 33 & 1.3141 & 0.8800 & 62 & [50, 71] & \citep{Rolfes2019} \\
2017--2018 & A(H3N2) & 166 & 1.6792 & 1.1411 & 22 & [12, 31] & \citep{Rolfes2019} \\
2017--2018 & B/Yamagata & 38 & 1.0612 & 0.7985 & 48 & --- & \citep{Rolfes2019} \\
2018--2019 & A(H1N1)pdm09 & 68 & 1.7398 & 1.0124 & 44 & [37, 51] & \citep{Chung2020} \\
2018--2019 & A(H3N2) & 47 & 1.9652 & 1.5111 & 9 & [0, 20] & \citep{Chung2020} \\
2019--2020 & A(H1N1)pdm09 & 47 & 1.8124 & 0.9147 & 30 & [21, 39] & \citep{Tenforde2021} \\
2019--2020 & B/Victoria & 108 & 1.5570 & 1.4553 & 45 & [37, 52] & \citep{Tenforde2021} \\
2021--2022 & A(H3N2) & 86 & 2.1416 & 1.4638 & 36 & [20, 49] & \citep{Price2023} \\
2022--2023 & A(H1N1)pdm09 & 39 & 1.8935 & 1.0199 & 47 & [1, 72] & \citep{Chung2024} \\
2023--2024 & A(H1N1)pdm09 & 54 & 1.2748 & 0.8464 & 29 & [15, 41] & \citep{Chung2025} \\
2023--2024 & A(H3N2) & 13 & 1.8180 & 1.4362 & 30 & [8, 47] & \citep{Chung2025} \\
2023--2024 & B/Victoria & 4 & 1.0127 & 0.5581 & 74 & [65, 81] & \citep{Chung2025} \\
2024--2025 & A(H1N1)pdm09 & 73 & 1.6211 & 1.1130 & 37 & [24, 48] & \citep{Chung2026} \\
2024--2025 & A(H3N2) & 66 & 1.4044 & 1.2941 & 27 & [14, 39] & \citep{Chung2026} \\
2024--2025 & B/Victoria & 3 & 1.1396 & 0.7173 & 40 & [12, 59] & \citep{Chung2026} \\
\bottomrule
\end{tabular}
}
\end{table}

\subsection{Vaccine strain records and GISAID sequence sets}
\label{sec:sequencerecords}

Tables~\ref{tab:targets} and \ref{tab:episets} identify the vaccine targets
and registered preseason sequence sets. The 23 distinct seed strains supply
43 season--component assignments. Five assignments use records deposited
after the corresponding season: B/Victoria in 2010--2011, 2011--2012,
2014--2015 and 2015--2016, and H3N2 in 2012--2013. Although strain names
match the recommendations, these later records limit strict claims about
historical data availability. Target mapping was checked manually.

\begin{table}[!htbp]
\centering
\caption{Vaccine target (seed strain) inventory. For each season and component
the target is the named WHO/NIBSC/NYMC/CDC reassortant seed strain where one
exists on GISAID, egg-passaged, rather than an arbitrary egg-passaged clinical
isolate, because the named reassortant is the sequence used to manufacture the
vaccine. Accessions and passage histories are as recorded in GISAID. The
mapping from each WHO February recommendation to a target sequence was verified
manually rather than programmatically.}
\label{tab:targets}
\scriptsize
\begingroup\sloppy
\begin{tabular}{@{}>{\raggedright\arraybackslash}p{0.255\textwidth}@{\hspace{6pt}}l@{\hspace{6pt}}l@{\hspace{6pt}}>{\raggedright\arraybackslash}p{0.30\textwidth}@{}}
\toprule
Vaccine seed strain & GISAID EPI\_ISL & Passage & Seasons and components \\
\midrule
A/Missouri/11/2025 (IVR-279) & \texttt{EPI\_ISL\_20215073} & Spf3E1D6 & \textbf{2026-27 prediction}, A(H1N1)pdm09 (egg) \\
A/Darwin/1454/2025 (IVR-283) & \texttt{EPI\_ISL\_20356066} & E2/D7 & \textbf{2026-27 prediction}, A(H3N2) (egg) \\
B/Tokyo/EIS13-175/2025 & \texttt{EPI\_ISL\_20373027} & E3 & \textbf{2026-27 prediction}, B/Victoria (no named reassortant yet) \\
A/California/7/2009 & \texttt{EPI\_ISL\_391380} & E2E6 & 2010-11 -- 2015-16, A(H1N1)pdm09 \\
B/Brisbane/60/2008 (NIB-58) & \texttt{EPI\_ISL\_232491} & E4 & 2010-11, 2011-12, 2013-14 -- 2017-18, B/Victoria \\
A/Perth/16/2009 (NIB-64) & \texttt{EPI\_ISL\_70443} & E4/E2 & 2010-11, 2011-12, A(H3N2) \\
A/Victoria/361/2011 (IVR-165) & \texttt{EPI\_ISL\_168948} & E6 & 2012-13, A(H3N2) \\
B/Wisconsin/01/2010 (BX-41A) & \texttt{EPI\_ISL\_83975} & EX & 2012-13, B/Yamagata \\
A/Texas/50/2012 (X-223A) & \texttt{EPI\_ISL\_132571} & E5E8 & 2013-14, 2014-15, A(H3N2) \\
B/Massachusetts/2/2012 (BX-51C) & \texttt{EPI\_ISL\_402396} & E3E7/E1 & 2013-14, 2014-15, B/Yamagata \\
B/Phuket/3073/2013 & \texttt{EPI\_ISL\_291275} & E6 & 2015-16 -- 2024-25, B/Yamagata (every season it's the vaccine strain) \\
A/Switzerland/9715293/2013 (X-247) & \texttt{EPI\_ISL\_169603} & E4E9+1 & 2015-16, A(H3N2) \\
A/Hong Kong/4801/2014 (X-263B) & \texttt{EPI\_ISL\_222828} & EX/E3 & 2016-17, 2017-18, A(H3N2) \\
A/Michigan/45/2015 (X-275B) & \texttt{EPI\_ISL\_233044} & EX & 2017-18, 2018-19, A(H1N1)pdm09 \\
B/Colorado/06/2017 & \texttt{EPI\_ISL\_276408} & E5 & 2018-19, 2019-20, B/Victoria (and combined B in 2018-19) \\
A/Singapore/INFIMH-16-0019/2016 (X-307) & \texttt{EPI\_ISL\_293533} & E5E2E9 & 2018-19, A(H3N2) \\
A/Brisbane/02/2018 (IVR-190) & \texttt{EPI\_ISL\_20212507} & E3/D8/E2 & 2019-20, A(H1N1)pdm09 \\
A/Victoria/2570/2019 (X-379A) & \texttt{EPI\_ISL\_14809235} & EX & 2021-22, 2022-23, A(H1N1)pdm09 \\
A/Cambodia/e0826360/2020 (X-361A) & \texttt{EPI\_ISL\_2378896} & EX & 2021-22, A(H3N2) \\
B/Austria/1359417/2021 (BVR-26) & \texttt{EPI\_ISL\_5933268} & Ex & 2022-23, 2023-24, 2024-25, B/Victoria \\
A/Darwin/9/2021 (NYMC X-369A) & \texttt{EPI\_ISL\_5522637} & EX & 2022-23, 2023-24, A(H3N2) (egg) \\
A/Victoria/4897/2022 (IVR-238) & \texttt{EPI\_ISL\_17086692} & E3/E6 & 2023-24, 2024-25, A(H1N1)pdm09 (egg) \\
A/Thailand/8/2022 (IVR-237) & \texttt{EPI\_ISL\_18372670} & E3/D7 & 2024-25, A(H3N2) (egg) \\
\bottomrule
\end{tabular}
\endgroup

\end{table}

\begin{table}[!htbp]
\centering
\caption{GISAID EPI\_SET accessions and DOIs for the preseason searches of each
season contributing to the cohort or the forecast. Each DOI resolves to the
acknowledgement table for that sequence set. The 2026--2027 entries are the
window-closed searches of 9 October 2026 used for the frozen forecast.}
\label{tab:episets}
\scriptsize
\begin{tabular}{llll}
\toprule
Season & Component & GISAID EPI\_SET & DOI \\
\midrule
2010--2011 & A(H1N1)pdm09 & \texttt{EPI\_SET\_260916yo} & \doilink{10.55876/gis8.260916yo} \\
 & A(H3N2) & \texttt{EPI\_SET\_260916yt} & \doilink{10.55876/gis8.260916yt} \\
 & B (both lineages) & \texttt{EPI\_SET\_260916vw} & \doilink{10.55876/gis8.260916vw} \\
2011--2012 & A(H1N1)pdm09 & \texttt{EPI\_SET\_260916er} & \doilink{10.55876/gis8.260916er} \\
 & A(H3N2) & \texttt{EPI\_SET\_260916mg} & \doilink{10.55876/gis8.260916mg} \\
 & B (both lineages) & \texttt{EPI\_SET\_260916kv} & \doilink{10.55876/gis8.260916kv} \\
2012--2013 & A(H1N1)pdm09 & \texttt{EPI\_SET\_260916sc} & \doilink{10.55876/gis8.260916sc} \\
 & A(H3N2) & \texttt{EPI\_SET\_260916fs} & \doilink{10.55876/gis8.260916fs} \\
 & B (both lineages) & \texttt{EPI\_SET\_260916zm} & \doilink{10.55876/gis8.260916zm} \\
2014--2015 & A(H3N2) & \texttt{EPI\_SET\_260916eb} & \doilink{10.55876/gis8.260916eb} \\
 & B (both lineages) & \texttt{EPI\_SET\_260916fh} & \doilink{10.55876/gis8.260916fh} \\
2015--2016 & A(H1N1)pdm09 & \texttt{EPI\_SET\_260916hw} & \doilink{10.55876/gis8.260916hw} \\
 & A(H3N2) & \texttt{EPI\_SET\_260916cv} & \doilink{10.55876/gis8.260916cv} \\
 & B (both lineages) & \texttt{EPI\_SET\_260916qy} & \doilink{10.55876/gis8.260916qy} \\
2016--2017 & A(H3N2) & \texttt{EPI\_SET\_260916ru} & \doilink{10.55876/gis8.260916ru} \\
 & B (both lineages) & \texttt{EPI\_SET\_260916az} & \doilink{10.55876/gis8.260916az} \\
2017--2018 & A(H1N1)pdm09 & \texttt{EPI\_SET\_260916rg} & \doilink{10.55876/gis8.260916rg} \\
 & A(H3N2) & \texttt{EPI\_SET\_260916pu} & \doilink{10.55876/gis8.260916pu} \\
 & B (both lineages) & \texttt{EPI\_SET\_260916ce} & \doilink{10.55876/gis8.260916ce} \\
2018--2019 & A(H1N1)pdm09 & \texttt{EPI\_SET\_260916nr} & \doilink{10.55876/gis8.260916nr} \\
 & A(H3N2) & \texttt{EPI\_SET\_260916ux} & \doilink{10.55876/gis8.260916ux} \\
 & B (both lineages) & \texttt{EPI\_SET\_260916nb} & \doilink{10.55876/gis8.260916nb} \\
2019--2020 & A(H1N1)pdm09 & \texttt{EPI\_SET\_260916es} & \doilink{10.55876/gis8.260916es} \\
 & B (both lineages) & \texttt{EPI\_SET\_260916ps} & \doilink{10.55876/gis8.260916ps} \\
2021--2022 & A(H3N2) & \texttt{EPI\_SET\_260916wu} & \doilink{10.55876/gis8.260916wu} \\
2022--2023 & A(H1N1)pdm09 & \texttt{EPI\_SET\_261014xz} & \doilink{10.55876/gis8.261014xz} \\
2023--2024 & A(H1N1)pdm09 & \texttt{EPI\_SET\_260916kg} & \doilink{10.55876/gis8.260916kg} \\
 & A(H3N2) & \texttt{EPI\_SET\_260916aq} & \doilink{10.55876/gis8.260916aq} \\
 & B (both lineages) & \texttt{EPI\_SET\_260916tr} & \doilink{10.55876/gis8.260916tr} \\
2024--2025 & A(H1N1)pdm09 & \texttt{EPI\_SET\_260916hs} & \doilink{10.55876/gis8.260916hs} \\
 & A(H3N2) & \texttt{EPI\_SET\_260916du} & \doilink{10.55876/gis8.260916du} \\
 & B (both lineages) & \texttt{EPI\_SET\_260916xc} & \doilink{10.55876/gis8.260916xc} \\
2026--2027 & A(H3N2) & \texttt{EPI\_SET\_261009ud} & \doilink{10.55876/gis8.261009ud} \\
 & A(H1N1)pdm09 & \texttt{EPI\_SET\_261009xz} & \doilink{10.55876/gis8.261009xz} \\
 & B/Victoria & \texttt{EPI\_SET\_261009tf} & \doilink{10.55876/gis8.261009tf} \\
\bottomrule
\end{tabular}

\end{table}

\subsection{Verification of published VE estimates and figure-derived values}
\label{sec:vechecks}

Each outcome was checked against its original article before model
refitting. Two of 28 assembled observations were excluded: the recorded
2017--2018 B/Victoria value of 76\% was not found as a lineage-specific
estimate in \citet{Rolfes2019}, and the recorded 2018--2019 B/Yamagata value
of 34\% represented combined influenza B rather than that lineage
\citep{Chung2020}.

Twenty-five retained estimates were printed numbers. These include the four
2015--2016 estimates printed in Figure 1 of \citet{Jackson2017}. The remaining
estimate, 48\% for 2017--2018 B/Yamagata, was read from the unlabelled point
in Figure 2 of \citet{Rolfes2019}. Applying the same reading procedure to
four points with printed values recovered 38.1\%, 22.0\%, 62.0\% and 49.9\%
against 38\%, 22\%, 62\% and 50\%, respectively. Omitting the figure-derived
estimate left 25 observations and gave the pooled topological fit
$\widehat{\VE}=77.63-23.54D_A$ ($R^2=0.301$; slope-permutation $p=0.0050$),
compared with $R^2=0.304$ and $p=0.0033$ in the primary analysis.
The pooled association therefore did not depend on retaining this single
figure-derived estimate.

Only point estimates entered the models. Confidence intervals were not
available for five observations; among the 21 recorded intervals, median
width was 22 percentage points (mean 26.9; range 14--71). The interval held
for 2015--2016 H3N2 came from an agency archive rather than the article and
is not shown as a primary-source interval in Table~\ref{tab:cohort}. The
2022--2023 estimate is late-season rather than a conventional full-season
estimate \citep{Chung2024}. Clinical estimation uncertainty was not propagated
through the regression.

\subsection{Word lengths, radius grids and checks of distance calculations}
\label{sec:implementation}

Both distances used HA proteins, the 20 amino acids and word lengths one to
three. The protein examples of \citet{HozumiWei2024} also use maximum word
length three. The Betti-curve implementation used only Betti-0
($b_k=0$ in that paper's Equation 16), with radii
$\{1,2,4,\ldots,N\}$, $N$ being the greater of the two sequence lengths. The
Grassmann implementation used 50 log-spaced radii on $[1,600]$, four channels
and rank 10. Both grids depart from the linear schedule $r_i=4^{\,k-1}i$ of the
source nucleotide implementation, for the reason given in
Section~\ref{sec:instantiation}; the upper radius exceeds the longest protein
in this dataset (585 amino acids).

The Betti-curve calculation was checked against the published formulas and
internal consistency tests; it was not numerically compared with the
nucleotide-only reference code. The Grassmann calculation used the reference
feature, normalisation, SVD and distance routines of \citet{WangWei2025},
with the amino-acid word enumeration and radius grid supplied separately.
With the reference nucleotide settings, feature blocks matched to machine
precision across all four channels and word lengths one to three. Removing
rows that are zero in every sequence preserves nonzero features and the
subspace distances, provided shared coordinates and the RMS denominator
are preserved.

\subsection{Permutation tests and intervals for future VE}
\label{sec:statdetails}

For the two primary pooled associations, we used two-sided slope tests
with 20,000 permutations of the VE labels (NumPy random seed 0). For each
permutation, we refitted the regression and compared its absolute slope with
the observed absolute slope; the $p$ value was the fraction of permuted
slopes at least as large. This gave $p=0.0033$ for the topological distance
and $p=0.0009$ for the geometric distance. The sequence-exclusion and
printed-number-only sensitivity tests used the same procedure.

We report 95\% OLS future-observation prediction intervals. These account
for fitted-regression uncertainty and residual variation under the linear
model, rather than uncertainty in the mean fitted VE alone.

\subsection{Training seasons, test observations and distance scales}
\label{sec:temporaldetails}

The evaluation retained the original common test set for comparison with
the historical-mean reference. A season was eligible only if every viral
group in that season had an earlier VE observation. The training set also
had to contain more observations than the number of previously observed
viral groups plus one, ensuring estimable residual variance under the
original comparison protocol. These rules were retained when reporting
only the pooled models.

The 19 test observations came from 2014--2015 (one), 2016--2017 (three),
2017--2018 (three), 2018--2019 (two), 2019--2020 (two), 2021--2022 (one),
2022--2023 (one), 2023--2024 (three) and 2024--2025 (three).
The first three observations, from 2010--2011, 2011--2012 and 2012--2013,
provided the initial training set. The 2015--2016 season was excluded as a
whole test season because three viral groups first appeared that year;
its four observations entered training for subsequent seasons.
Thus, seven of the 26 observations contributed to training without being
scored as test predictions.

At each test season, we refitted both the regression and distance scales
using all eligible observations from earlier seasons. The Betti-curve scales
were estimated from their virus--vaccine comparisons. Grassmann scales
were estimated from the unique circulating and vaccine protein sequences
in the same training data. These scales were then applied to both training
and test distances. No sequence from the test season or a later season
contributed to scale estimation.

The reference predicted the earlier mean VE for the relevant subtype or
lineage. Errors were calculated once per season--component observation,
with equal weight across the test observations. MAE averages absolute
errors; RMSE takes the square root of the mean squared error. Coverage
counts observed VE estimates inside their prediction intervals, and mean
width averages the differences between interval limits. The relative MAE
reduction is $100(1-\mathrm{MAE}_{\mathrm{model}}/\mathrm{MAE}_{\mathrm{reference}})$.

For Grassmann scale estimation, the RMS denominator counted all $20^k$
possible amino-acid words, including zero rows, so later sequences could
not alter the row-count convention. A check using the original stored-row
denominator gave a geometric MAE reduction of 6.2\%, compared with 5.2\%
under the full-alphabet convention. Both calculations used the same 19
test observations; this check concerns preprocessing. Neither reduction is
distinguishable from zero (Section~\ref{sec:pooledtemporal}), so the
convention does not change the conclusion.
The training-only evaluation was implemented retrospectively.
Original pooled fits and frozen forecasts retain the project-wide scales
described in Section~\ref{sec:instantiation}.

The split follows the seasons in which specimens and clinical outcomes
were collected. It does not reconstruct exactly which sequence records or
published VE estimates were publicly available at each historical forecast
date. Consequently, this is a retrospective evaluation with season-ordered
training and testing; prospective performance remains to be established.

\subsection{Matching sampling windows and omitting individual pairs or seasons}
\label{sec:windowdetails}
\label{sec:discussiondetails}

The topological window comparison used 18 matched observations across eight
seasons, with the same vaccine targets and scales in both windows. Scales
were estimated from the preseason sequences in the combined 30,464-record
collection and applied to both periods. This collection combines the original
paired-window search with a later recovery of A(H1N1)pdm09 records, after
removing the records flagged by the nonseasonal-virus checks. Correlations were calculated before
and after subtracting group means separately within each window; uncertainty
was assessed by resampling whole seasons. This comparison does not use VE
labels. Enrollment-period sequences did not enter the primary VE regressions.

The component-demeaned correlation ranged from 0.755 to 0.945 after deleting
one pair and from 0.783 to 0.948 after deleting one season. The correlation
remained positive and high in each deletion, so no single pair or season
was solely responsible for the overall relationship. Within-group
correlations were 0.98 for H1N1 ($n=5$), 0.93 for B/Victoria ($n=4$), 0.70
for H3N2 ($n=6$) and $-0.30$ for B/Yamagata ($n=3$). In 2018--2019 H3N2,
the preseason distance was within 0.06 of its group mean, while the
enrollment-period distance was 0.46 above its group mean. Thus, high overall
correlation did not imply agreement in every group or season.

The enrollment-period collection contains 29,461 sequences. No GISAID
EPI\_SET identifier was registered when these searches were performed, so these
sequences are not covered by the identifiers in Supplementary
Table~\ref{tab:episets}. Their accessions are instead listed by season, viral
group and sampling window in \texttt{metric\_validation/paired\_window\_accessions.csv}
in the repository, from which the corresponding GISAID records and their
contributor acknowledgements can be retrieved. These searches nonetheless
remain less fully documented than the preseason datasets. The matched
seasons were selected by data availability rather than random sampling.

\end{document}